\documentclass[leqno,11pt]{amsart}
\usepackage{graphicx}
\usepackage{amscd}
\usepackage{amsmath}
\usepackage{caption}
\usepackage{amsfonts}
\usepackage{amssymb}
\usepackage{mathrsfs}
\usepackage{multicol}
\usepackage{caption}
\usepackage{subcaption}
\usepackage{color}
\numberwithin{equation}{section}

\newcommand{\be}{\begin{equation}}
\newcommand{\ee}{\end{equation}}
\newcommand{\ben}{\begin{eqnarray*}}
\newcommand{\een}{\end{eqnarray*}}

\newtheorem{theorem}{Theorem}
\newtheorem{lemma}{Lemma}
\newtheorem{remark}{Remark}
\newtheorem{corollary}{Corollary}
\newtheorem{proposition}{Proposition}
\newtheorem{definition}{Definition}
\allowdisplaybreaks
\definecolor{darkgreen}{rgb}{0.09, 0.45, 0.27}
\definecolor{debianred}{rgb}{0.84, 0.04, 0.33}
\definecolor{orange}{rgb}{0.98, 0.6, 0.01}
\allowdisplaybreaks

\begin{document}


 \title[Generalized quasilinear Schr\"odinger equations in $\mathbb R^N$]
{Multiplicity results for generalized quasilinear critical Schr\"odinger equations in $\mathbb R^N$}

\author[L. Baldelli]{Laura Baldelli} \email{lbaldelli@impan.pl}

\author[R. Filippucci]{Roberta Filippucci} \email{roberta.filippucci@unipg.it}

\address[Filippucci]{Department of Mathematics -- University of Perugia --
Via Vanvitelli 1 -- 06123 Perugia,  Italy}
\address[Baldelli]{Institute of Mathematics -- Polish Academy of Sciences --
ul. Sniadeckich 8 -- 00-656 Warsaw, Poland}

\begin{abstract}
Multiplicity results are proved for solutions both with positive and negative energy, as well as nonexistence results, of a generalized quasilinear Schr\"odinger potential free equation in the entire $\mathbb R^N$  involving a nonlinearity which combines a power-type term at a critical level with a subcritical term, both with weights. 
The equation has been derived from models of several physical phenomena such as superfluid film in plasma physics as well as the self-channelling of a high-power ultra-short laser in matter.

Proof techniques, also in the symmetric setting, are based on variational tools, including concentration compactness principles, to overcome lack of compactness, and the use of a change of variable in order to deal with a well defined functional.
\end{abstract}

\keywords
{variational methods, generalized quasilinear Schr\"odinger equations, multiplicity results, concentration compactness\\
\phantom{aa} 2020 AMS Subject Classification: Primary: 35J62
Secondary: 35J60; 35J20.}

\maketitle

\section{Introduction}\label{intro}

In this paper, we are interested in multiplicity results for nontrivial weak solutions, both with negative and positive energy of the following generalized quasilinear Schr\"odinger equation involving a critical term
\begin{equation}\label{Psc}
-\Delta_{p}u-\frac{\alpha}{2}\Delta_p(|u|^\alpha)|u|^{\alpha-2}u=\lambda V(x)|u|^{k-2}u
+\beta K(x)|u|^{\alpha p^{*}-2}u \quad\text{in }\mathbb{R}^N,
\end{equation}
where $\Delta_{m}u=\mbox{div}(|Du|^{m-2}D u)$ is the $m$-Laplacian of $u$, $\alpha>1$, $\beta, \lambda>0$, $N\ge 3$,
 $1<p<N$, $p^{*}=Np/(N-p)$ is the critical Sobolev's exponent,
 the exponent $k$ is such that $\alpha<k<\alpha p^*$ and the weights are nontrivial and satisfy
\begin{equation}\label{V(x)1}
0\le V\in  L^{r}(\mathbb{R}^{N})\cap C(\mathbb R^N), \quad  r=\frac{\alpha p^{*}}{\alpha p^{*}-k},
\end{equation}
\begin{equation}\label{K(x)1}
K\in  L^{\infty}(\mathbb{R}^{N})\cap C(\mathbb R^N).
\end{equation}
In particular, we emphasize that the weight $K$ can change sign.

Actually, in the last section of the paper, we deal also with the singular case, given by $0<\alpha<1$, and, in particular, thanks to a technique used to attach the case $\alpha>1$,  we succeed in extending some previous results proved in \cite{BFsc}. Our multiplicity results are a first contribution in the study of critical generalized quasilinear Schr\"odinger equations involving $p$-Laplacian and general $\alpha>1$. To the best of our knowledge, they  are new even in the Laplacian case.

Solutions of problems of type \eqref{Psc} are related to the existence of standing wave solutions for Schr\"odinger  equations of the form
\begin{equation}\label{Sch}
i\partial_t \psi = -\Delta \psi+ W(x)\psi- \varphi(|\psi|^2)\psi-\kappa \Delta \varrho(|\psi|^2)\varrho'(|\psi|^2)\psi,
\end{equation}
where $\psi:\mathbb R\times \mathbb R^N \to \mathbb C$, $ W$ is a given potential, $\kappa\in\mathbb R$ and $\varphi,\varrho$ real functions of mainly pure power form. 
By using the well known energy methods, the semilinear subcase  $\kappa=0$ of \eqref{Sch} was studied deeply in \cite{JT}, see also and the references therein.  

Equation \eqref{Sch}, according to  different types of $\varrho$, has been derived from models of several physical phenomena. For instance, when $\varrho(s)=s$, then \eqref{Sch} is interpreted as
the superfluid film equation in plasma physics, while if  $\varrho(s) = (1+s)^{1/2}$, equation \eqref{Sch} models the self-channelling of a high-power ultra short laser in matter. We refer to \cite{LWW3}, for applications, in the theory of Heisenberg ferromagnets and magnons and to \cite{barjean}, where  models for binary mixtures of ultracold quantum gases are treated. 

It is worth pointing out that there are many difficulties in treating this class of generalized critical quasilinear Schr\"odinger equations \eqref{Psc} in $\mathbb R^N$, such as lack of compactness, the presence of the term  $\Delta_p(|u|^\alpha)|u|^{\alpha-2}u$ which prevents us from working directly in a classical working space, to manage weights and finally to deal with an "apparently" supercritical term $K(x)|u|^{\alpha p^{*}-2}u$. In light of this, challenging tasks, tricky to be managed, appear.

Actually, when $\alpha>1$, the exponent $\alpha p^{*}$ plays the role of the critical exponent, differently from the case $0<\alpha<1$ in which $p^*$ behaves like the critical exponent, as it appears in the critical equation \eqref{Pscc} below.
In particular,  Section \ref{prel3} is aimed to justify the appearance of the new critical  exponent $\alpha p^{*}$ in terms of the compactness of a suitable  embedding and of the nonexistence of solutions beyond $\alpha p^*$. In the spirit of Theorem 2.11 in \cite{AW2012}, where the semilinear version ($p=2$) of  \eqref{Psc} is investigated.
We refer also to \cite{DPY16}, \cite{DS} and \cite{WLA}.

 In addition,  the choice of the parameter $\alpha$ gives rise to a different nature of the equation under consideration since for $\alpha>1$,  the term $\Delta_p(|u|^\alpha)|u|^{\alpha-2}u$, $p>1$, is degenerate at $u=0$, while for $0<\alpha<1$ it becomes singular when $u=0$.

In literature, critical Schr\"odinger equations in $\mathbb R^N$ are mainly studied in their physical relevance, corresponding to the situation $p=2$ and when a potential term  is involved.
The degenerate case, $\alpha>1$, is mostly studied when $\alpha=2$, see \cite{SV10}, 
where existence results can be deduced by applying the Mountain Pass Theorem in the superlinear subcritical case $4 (=\alpha p)<k<2\cdot 2^*$ of \eqref{Psc}. 
As for equations  where the reaction combines the multiple effects generated by a singular term and a critical term  even with nonhomogeneous operators and in bounded domains,  we refer to  the recent works  \cite{Kumar} and  \cite{prrsing}.  While multiplicity results for fractional magnetic nonlinear Schr\"odinger equations, both in the critical case and in the supercritical case, are considered in the recent paper \cite{ambr}.

For existence results in the  critical semilinear case when $\alpha>1$ in \eqref{Psc} very few is known in the entire interval $\alpha<k<\alpha p^*$. Nonlinearities only subcritical are studied  in \cite{ASW14} when $2^*<k<2^*\alpha$, while \cite{LW02} deals with the case $2\le k<2^*\alpha$.
In the critical case, see also \cite{SW16}, \cite{DPY16} where a general nonlinearity is included.
For quasilinear critical Schr\" odinger equations with $\alpha=2$, a potential and involving a Choquard term, we refer to \cite{LiZh}.

Passing to the general quasilinear case, $1<p<N$,  with  $\alpha>1$, 
it seems that  there are no multiplicity results for quasilinear critical Schr\"odinger problems on $\mathbb R^N$ when a potential term is not involved in the equation  and the growth rate of the subcritical term is less than $\alpha p$.  
Motivated by this observation and with the sake of  completing the picture started in \cite{BFsc}, where the singular case $0<\alpha<1$ in \eqref{Psc} is treated, in this paper we investigate the case $\alpha<k<p$ in the following theorem, where $E_\lambda$ stands for the standard energy functional associated to equation \eqref{Psc}, see Section~\ref{prel3}.

\begin{theorem}\label{1.1}
Assume $N\ge3$, $1<\alpha<k< p$. Let $V$, $K$ satisfying \eqref{V(x)1} and \eqref{K(x)1}, respectively. Then, 
\begin{itemize}
\item[(i)] For any $\lambda>0$, there exists $\beta^*>0$ such that for any $0<\beta<\beta^*$, then equation \eqref{Psc} has infinitely many nontrivial solutions $(u_n)_n\subset D^{1,p}(\mathbb R^N)$  such that $E_\lambda(u_n)<0$ and $\|u_n\|\to0$ as $n\to\infty$.
\item[(ii)]For any $\beta>0$, there exists $\lambda^*>0$ such that for any $0<\lambda<\lambda^*$, then equation \eqref{Psc} has infinitely many nontrivial solutions $(u_n)_n\subset D^{1,p}(\mathbb R^N)$  such that $E_\lambda(u_n)<0$ and $\|u_n\|\to0$ as $n\to\infty$.
\end{itemize}
\end{theorem}

The proof of the above multiplicity result relies on the concentration compactness principle, the truncation of the energy functional and the theory of Krasnosel'skii genus. In particular, we cannot manage directly the energy functional $E_\lambda$ associated with \eqref{PS} since it might be not well defined, but we need to perform a suitable change of variables, in the spirit of \cite{CJ}, to achieve a nice functional. This causes some further obstacles in recovering in some sense compactness beyond $p^*$, since, as noted above, $\alpha p^*>p^*$.
While the behaviour of the norm of the solutions follows from the application of the symmetric Mountain Pass Theorem, see \cite{AR, C, K5, R}.


Moreover, we also need to restrict the entire sublinear interval $\alpha<k<\alpha p$ to $\alpha<k<p$ due to the shape of the transformed energy functional which yields tricky estimates, cfr. Remark \ref{k_restr}. As far as we know,  in \cite{WWL9} there is an attempt to cover the entire interval $\alpha<k<\alpha p$ when $K\equiv1$, by using an inequality which seems difficult to be verified.

In the second part of our paper, we study equation \eqref{Psc}  in a symmetric setting related to a subgroup $T$ of the group $O(N)$ of orthogonal
 linear transformations in $\mathbb{R}^N$. In turn, it is possible to define 
$$|T|:=\inf_{x\in\mathbb{R}^N, \,\ x\ne0}|T_{x}|,$$
where $|T_x|$ the cardinality of a $T$-orbit $T_x:= \{\tau x | \tau \in T\}$ with
$|T_\infty|:=1$ (so that $|T_0|=1$).

We make use of $T$-symmetric functions $f:\Omega\to \mathbb R^N$, i.e. $f(\tau x)=f(x)$ for all
$\tau \in T$ and $x\in \Omega$, with $\Omega$ open $T$-symmetric subset of $\mathbb R^N$ 
(i.e. if $x\in\Omega$, then $\tau x\in\Omega$ for all $\tau\in T$). For example, even functions are $T$-symmetric functions with $T=\{id,-id\}$,  thus $|T|=2$, and radially symmetric functions
are $T$-symmetric functions with $T=O(N)$, thus $|T|=\infty$.
Then,  we denote with $D^{1,p}_{T}(\mathbb R^N)$ the subspace of $D^{1,p}(\mathbb R^N)$ consisting
of all $T$-symmetric functions, cfr. Section \ref{prel3}.

A pioneering paper about critical symmetric problems in the entire $\mathbb R^N$ is \cite{BCS} by Bianchi, Chabrowskii and Szulkin, where existence and multiplicity results are obtained for $p=2$ under appropriate assumptions on the single weight involved and on the group $T$, no subcritical terms are included.
In this context, we mention also \cite{cc} where  the symmetry of solutions is shown by applying the improved "moving plane" method.

Differently from Theorem \ref{1.1}, where no assumptions on the sign for the weight $K$ are required, in the next result, nonnegativity for the weight $K$ is needed.
Furthermore, we consider solutions with positive energy, with no additional restrictions on $\alpha, p, k$, indeed the range for $\alpha$ and $k$ is the largest possible.

\begin{theorem}\label{T_E>0}
Assume $N\ge 3$, $1<\alpha<k<\alpha p^*$. Let $V$, $K$ be $T$-symmetric functions satisfying \eqref{V(x)1}, \eqref{K(x)1}, with $K$ nonnegative (nontrivial) in $\mathbb R^N$. If
\begin{equation}\label{K0Kinfty=0}
K(0)=K(\infty)=0\quad \text{and}\quad |T|=\infty,
\end{equation}
where $K(\infty)=\limsup_{|x|\to\infty} K(x)$. 
Then, for all $\lambda, \beta>0$ equation \eqref{Psc} possesses
infinitely many solutions $(u_n)_n\subset D^{1,p}_T(\mathbb R^N)$ with positive energy such that $E_\lambda(u_n)\to\infty$ as $n\to\infty$.
\end{theorem}

The main ingredient used in the proof of Theorem \ref{T_E>0} is the Fountain Theorem, which requires the Palais Smale property for the functional at any positive level. This is in force by virtue of the crucial assumption \eqref{K0Kinfty=0}.

The paper is structured as follows.
The description of the functional setting, the reformulation of the problem by a suitable change of variable and motivations on the critical exponent $\alpha p^*$ are contained in Section \ref{prel3}. 
Section \ref{psss} encloses compactness properties thanks to which we overcome the lack of compactness, while in Section \ref{truncated} we perform a deep analysis on the possible behaviours of the energy functional and of its truncated version, this latter introduced to restore the boundedness from below.
The proof of Theorem \ref{1.1} is disclosed in Section \ref{main1}, whereas in Section \ref{main2} we deal with solutions with positive energy and we develop the proof  Theorem \ref{T_E>0} in the symmetric setting described before.
Finally, Section \ref{sing}, by using some ideas given in Section \ref{psss}, is devoted to extending Theorem 1.1 in \cite{BFsc}, relative to the singular case $0<\alpha<1$, covering the entire interval $1<k<\alpha p$.

\section{Preliminaries}\label{prel3}
In this section, we introduce the main notations and we present preliminary results useful for the proof of the main
theorems of the paper, given in Sections \ref{main1}, \ref{main2}.

Let $D^{1,p}(\mathbb R^N)$ be the closure of $C_0^\infty(\mathbb R^N)$ with respect to the norm
$\|u\|_{D^{1,p}}=\|Du\|_p$, where $\|\cdot\|_{p}$ is the $L^p$ norm in $\mathbb R^N$. In particular, we can define the reflexive Banach space $D^{1,p}(\mathbb{R}^N)=\{u\in L^{p^{*}}(\mathbb{R}^N): Du\in L^{p}(\mathbb{R}^N)\}$.

The Euler Lagrange functional associated with equation \eqref{Psc} is the following
\begin{equation}\label{Sf}
E_{\lambda}(u)=\frac{1}{p}\int_{\mathbb R^N}g(u)^p|D u|^{p} dx
-\frac{\lambda}{k}\int_{\mathbb{R}^N}V|u|^{k}dx-\frac{\beta}{\alpha p^{*}}\int_{\mathbb{R}^N}K|u|^{\alpha p^{*}}dx, 
\end{equation}
for $u\in D^{1,p}(\mathbb{R}^N)$, where
\begin{equation}\label{gpicc}
g(t)=\Biggl[1+\frac{\alpha^p}{2}|t|^{p(\alpha-1)} \Biggr]^{1/p}, \quad t\in \mathbb R.
\end{equation}
Due to the appearance of the coercive term $g$, indeed $g(t)\to\infty$ as $|t|\to\infty$, when $\alpha>1$, 
the functional $E_{\lambda}$ may be not well defined in $D^{1,p}(\mathbb R^N)$, so we cannot apply variational methods 
or nonsmooth critical point theory to deal directly with \eqref{Sf}. 
For example, from \cite{S08}, if we consider
$u_e(x)= |x|^{-(N-p)/2p}$, $x\in B_1\setminus\{0\}$
then $u_e\in D^{1,p}(\mathbb R^N)$ but $g(u_e)^p|Du_e|^p \notin L^1({\mathbb R^N})$ for any $\alpha\ge2$.
To overcome this difficulty, we make a change of variables developed in \cite{CJ}, following an idea in \cite{LWW3}, precisely
\begin{equation}\label{ggra}
v=G(u)=\int_0^u g(z) dz.
\end{equation}
In particular, the function $g$ defined in \eqref{gpicc} is an even function in $\mathbb R$, $g(0)=1$, $g$ is increasing in $\mathbb R^+$ and decreasing in $\mathbb R^-$.
For any $\alpha>0$ and $t\in \mathbb R$, we have
\begin{equation}\label{G1a}|G(t)|\le \int_0^{|t|} \biggl( 1+\frac{\alpha}{2^{1/p}} y^{\alpha-1} \biggr) dy =|t|+\frac{1}{2^{1/p}}|t|^\alpha
\end{equation}
Thus, $G$ is well defined and continuous in $\mathbb R$. Moreover, $G$ is a strictly increasing function being $g\ge1$, $G(0)=0$ and $\lim_{|t|\to \infty}G(t)=\infty$. So, we can define $G^{-1}$, an invertible, odd and $C^1$ function such that $G^{-1}(s)\ge 0$ for any $s\ge0$.
Thanks to the change of variables described above, the energy functional $E_{\lambda}$ can be written by the following functional
\begin{equation}\label{F}
F_{\lambda}(v):=\frac{1}{p}\int_{\mathbb R^N}|D v|^{p} dx-\frac{\lambda}{k}
\int_{\mathbb{R}^N}V|G^{-1}(v)|^{k}dx-\frac{\beta}{\alpha p^{*}}\int_{\mathbb{R}^N}K|G^{-1}(v)|^{\alpha p^{*}}dx,
\end{equation}
for $v\in D^{1,p}(\mathbb{R}^N)$.
The proof of the regularity of $F_\lambda$ takes the following steps, starting with the properties of $g$ and $G$.
Especially, the following lemma holds.

\begin{lemma}\label{gprop}
Let $\alpha>1$. Then, it holds
\begin{enumerate}
\item[{\bf a)}]$\lim_{s\to 0} \frac{|G^{-1}(s)|}{|s|}=1$ thus
$\lim_{s\to0} \frac{|G^{-1}(s)|^\alpha}{|s|}=0$;
\item[{\bf b)}] $\lim_{s\to \infty} \frac{|G^{-1}(s)|^\alpha}{|s|}=2^{1/p}$ thus $\lim_{s\to \infty} \frac{|G^{-1}(s)|}{|s|}=0$;
\item[{\bf c)}] $|G^{-1}(s)|\le |s|$, for every $s\in\mathbb R$;
\item[{\bf d)}] $0\le\frac{g'(t) t}{g(t)}<\alpha-1$, for every $t\in\mathbb R$;
\item[{\bf e)}] $|G^{-1}(s)|\le |G^{-1}(s)g(G^{-1}(s))|< \alpha|s|$, for every $s\in\mathbb R$;
\item[{\bf f)}] Take $\wp\ge0$. Then the following hold for every $t\in\mathbb R$.

If $\wp>\alpha$ then $\frac{\wp-\alpha}{g(t)}< \frac{(\wp-1)g(t)-tg'(t)}{g^2(t)}\le \wp-1$.

If $1<\wp<\alpha$ then
$\wp-\alpha< \frac{(\wp-1)g(t)-tg'(t)}{g^2(t)}\le  \wp-1$.

\item[{\bf g)}]  $|G^{-1}(s)|^{\alpha}\le 2^{p-1}|s|$, for every $s\in\mathbb R$;
\item[{\bf h)}] $|G^{-1}(s)|^{\alpha}\ge |G^{-1}(1)|^{\alpha}|s|$, for every $s\in\mathbb R$ such that $|s|\ge 1$;
\item[{\bf i)}] $|G^{-1}(s)|^{\alpha-1}(G^{-1}(s))'\le C$ with $C>0$.
\end{enumerate}
\end{lemma}

\begin{proof}
Since $G^{-1}$ is odd, we only consider the case $s\ge 0$. Property  {\bf a)} is trivial by Hospital's rule, being $t=G^{-1}(s)$ so that $G'(t)=g(t)$ with $g(0)=1$.
Also {\bf b)} follows immediately again from Hospital's rule since 
$$\lim_{t\to \infty} \frac{t^\alpha}{G(t)}=\alpha\lim_{t\to \infty}\frac{t^{\alpha-1}}{g(t)}=2^{1/p},\qquad\lim_{t\to \infty} \frac{t}{G(t)}=\lim_{t\to \infty}\frac{1}{g(t)}=0,$$
where we have used also that 
\begin{equation}\label{org_g_0}
g(t)\sim \frac{\alpha}{2^{1/p}}t^{\alpha-1},\quad \text{ as } t\to \infty.\end{equation}
Condition {\bf c)} follows since $g(t)>1$ for all $t> 0$, thus $G(t)=\int_0^t g(z) dz\ge t$ for every $t\ge 0$.
To prove {\bf d)}, since $g$ is increasing in $\mathbb R^+$ and positive, we have
$$0\le\frac{g'(t) t}{g(t)}=(\alpha-1)\cdot\frac{\alpha^p t^{p(\alpha-1)}}
{2+\alpha^p t^{p(\alpha-1)}}< \alpha-1,$$
being $\alpha>1$.
To get {\bf e)}, multiply  {\bf d)} by $g(t)>0$ and integrate so that
$$(\alpha-1) G(t)>\int_0^{t} \bigl\{[g(z) z]'-g(z)\bigr\} dz=g(t) t- G(t),\quad t>0.$$
In turn, inequality {\bf e)} follows taking $s=G(t)$. 
For the proof of {\bf f)}, it is enough to multiply {\bf d)} by $-1/g(t)$ and then add  $(\wp-1)/g(t)$ so that
$$\frac{\wp-\alpha}{g(t)}< \frac{(\wp-1)g(t)-tg'(t)}{g^2(t)}\le \frac{\wp-1}{g(t)},$$
yielding {\bf f)} depending on the case.
In order to prove {\bf g)} and {\bf h)} take
$$\biggl(\frac{|G^{-1}(s)|^{\alpha}}{s}\biggr)'=\frac{G^{-1}(s)^{\alpha-1}[\alpha s-g(G^{-1}(s))G^{-1}(s)]}{s^2g(G^{-1}(s))}>0, \quad \text{for} \quad s\ge0$$
from {\bf e)}. Consequently, $|G^{-1}(s)|^{\alpha}/s$ is strictly increasing and  by {\bf b)} its limit at infinity is $2^{1/p}$, thus {\bf g)} follows immediately. In addition, we get $|G^{-1}(s)|^{\alpha}>|G^{-1}(1)|^{\alpha}s$ for $s\ge 1$ and {\bf h)} holds by virtue of symmetry. 
Finally, {\bf i)} follows from $G^{-1}(0)=0$ and \eqref{org_g_0}.
\end{proof}

\begin{remark}\label{gd1p} 
Note that, by definition, we have
\begin{equation}\label{eq1bb}
[G^{-1}(s)g(G^{-1}(s))]'=1+\frac{G^{-1}(s)g'(G^{-1}(s))}{g(G^{-1}(s))}\in (1,\alpha)
\end{equation}
thanks to Lemma \ref{gprop}-{\bf d)}.
Moreover, for any $v\in D^{1,p}(\mathbb R^N)$
$$
D[G^{-1}(v)g(G^{-1}(v))]=  \Biggl[1+\frac{G^{-1}(v)g'(G^{-1}(v))}{g(G^{-1}(v))}\Biggr] Dv$$
so, by using \eqref{eq1bb}, we obtain
\begin{equation}\label{eq1}
 |Dv|\le |D[G^{-1}(v)g(G^{-1}(v))]|\le\alpha |Dv|.
\end{equation}
In addition, since $g\ge 1$, then we have
$|DG^{-1}(v)|=|Dv|/g(G^{-1}(v))\le |Dv|$.
Thus, for any $v\in D^{1,p}(\mathbb R^N)$ we have $G^{-1}(v)g(G^{-1}(v)), G^{-1}(v)\in D^{1,p}(\mathbb R^N)$. 
On the other hand, if we take $v\in D^{1,p}(\mathbb R^N)$, by using the definition of $g$ in \eqref{gpicc}, we have
\begin{equation}\label{alphap}\begin{aligned}
|D&(|G^{-1}(v)|^{\alpha})|^p=\alpha^p |G^{-1}(v)|^{p(\alpha-1)} |DG^{-1}(v)|^p=\\
&=\alpha^p |G^{-1}(v)|^{p(\alpha-1)} \biggl|\frac{Dv}{g(G^{-1}(v))}\biggr|^p= 2\frac{\dfrac{\alpha^p}{2}|G^{-1}(v)|^{p(\alpha-1)}}{1+\dfrac{\alpha^p}{2}|G^{-1}(v)|^{p(\alpha-1)}}|Dv|^p\le 2|Dv|^p,
\end{aligned}\end{equation}
so that for any $v\in D^{1,p}(\mathbb R^N)$ also $|G^{-1}(v)|^\alpha\in D^{1,p}(\mathbb R^N)$.
\end{remark}

In what follows we make use of the next crucial lemma, for its proof we refer to Lemma 2.2 in \cite{BFsc} where we consider $|G^{-1}(\cdot)|^\alpha$ in place of $G^{-1}(\cdot)$ and, by Remark \ref{gd1p}, $|G^{-1}(v)|^\alpha\in D^{1,p}(\mathbb R^N)$ for every $v\in D^{1,p}(\mathbb R^N)$.

\begin{lemma}\label{g-1N}
Assume $v_n \rightharpoonup v$ in $D^{1,p}(\mathbb R^N)$, then $|G^{-1}(v_n)|^\alpha \rightharpoonup |G^{-1}(v)|^\alpha$ in $D^{1,p}(\mathbb R^N)$.
\end{lemma}

Actually, as described in the Introduction,  equation \eqref{Psc} is critical since the corresponding critical exponent in the nonlinearity is $\alpha p^*$, as soon as  $\alpha>1$. For the Laplacian case $p=2$, we refer to \cite{AW2012} where a detailed discussion in this direction is conducted, see also \cite{WLA}, \cite{DPY16}.
In order to justify this assumption, we first prove nonexistence beyond $\alpha p^*$, and then we recover continuity and local compactness until $\alpha p^*$.

In the following theorem, by the celebrated variational identity by Pucci and Serrin in \cite{PS_vi}, we immediately get that for $k\ge \alpha p^*$ nonexistence follows at all, generalizing the result in \cite{LWW4} where the authors studied the case $p=\alpha=2$ provided nonexistence of solutions in $H^1(\mathbb R^N)$ with $|Du|^2|u|^{2}\in L^1(\mathbb R^N)$.

\begin{theorem}
Equation \eqref{Psc} does not admit any solutions if the following hold
\begin{equation}\label{non_ex}
k\ge \alpha p^*, \quad x\cdot DV(x)\le0, \quad x\cdot DK(x)\le0
\end{equation}
\end{theorem}
\begin{proof}
By using the following Pucci and Serrin variational identity,  \cite{PS_vi},  in $\mathbb R^N$
$$\begin{aligned}0=\int_{\mathbb R^N}\biggl\{\mathcal{F}(x,u,Du)& \mbox{div}h+h_i\mathcal{F}_{x_i}(x,u,Du)-\biggl[\frac{\partial u}{\partial x_j} \frac{\partial h_j}{\partial x_i}+u \frac{\partial a}{\partial x_i}\biggr]\mathcal{F}_{p_i}(x,u,Du)\\ &-a \biggl[\frac{\partial u}{\partial x_i}\mathcal{F}_{p_i}(x,u,Du)+u\mathcal{F}_{u}(x,u,Du)\biggr] \biggr\}dx, \end{aligned}$$
for
$$\mathcal F(x,u,Du)=\frac 1p g(u)^p|Du|^p-\frac \lambda k V(x)|u|^k-\frac{\beta}{\alpha p^*}K(x)|u|^{\alpha p^*},$$
we get 
$$\begin{aligned}
&0=\int_{\mathbb R^N}\biggl\{\biggl[\frac 1p g(u)^p|Du|^p-\frac \lambda k V|u|^k-\frac{\beta}{\alpha p^*}K|u|^{\alpha p^*}\biggr]\,\mbox{div}h\\
&+h_i \biggl(-\frac{\lambda}{k}\frac{\partial V}{\partial x_i}|u|^k-\frac{\beta}{\alpha p^*} \frac{\partial K}{\partial x_i}|u|^{\alpha p^*} \biggr)-g(u)^p|Du|^{p-2}
\frac{\partial u}{\partial x_i}\frac{\partial u}{\partial x_j} \frac{\partial h_j}{\partial x_i}\\
&-g(u)^p|Du|^{p-2}u
\frac{\partial u}{\partial x_i}\frac{\partial a}{\partial x_i}
-a g(u)^{p-1}|Du|^{p}\bigl[g(u)+ug'(u)\bigr]+a \bigl [ \lambda V|u|^k
+\beta K|u|^{\alpha p^*}\bigr]\biggr\}dx.
\end{aligned}$$

Choosing 
$h(x)=x$ and $a(x)=N/\alpha p^*$
we arrive to
$$\begin{aligned}
0=&\int_{\mathbb R^N}\biggl\{\biggr[\frac Np g(u)^p|Du|^p-\frac{N\lambda}{k} V|u|^k-\frac{\beta N}{\alpha p^*}K|u|^{\alpha p^*}\biggr]\\
&-\frac{\lambda}{k}\frac{\partial V}{\partial x_i}|u|^k-\frac{\beta}{\alpha p^*} \frac{\partial K}{\partial x_i}|u|^{\alpha p^*}
-g(u)^p|Du|^{p-2}
\frac{\partial u}{\partial x_i}\frac{\partial u}{\partial x_j} \delta_{ij}\\
&-\frac{N}{\alpha p^*} g(u)^{p-1}|Du|^{p}\bigl[g(u)+ug'(u)\bigr]+\frac{\lambda N}{\alpha p^*} V|u|^k
+\frac{\beta N}{\alpha p^*} K|u|^{\alpha p^*}\biggr\}dx,
\end{aligned}$$
namely
\begin{equation}\label{VI_final}
\begin{aligned}
\int_{\mathbb R^N}\biggr\{
\frac{N-p}{\alpha p}
\bigl[(\alpha-1)g(u)-&ug'(u)\bigr] g(u)^{p-1}|Du|^p- \lambda N
\biggl(\frac{1}{k}-\frac{1}{\alpha p^*}\biggl)V|u|^k
\biggr\}dx
\\&=\int_{\mathbb R^N}\biggl\{
\frac{\lambda}{k}|u|^k\, x\cdot DV(x)\, +\frac{\beta}{\alpha p^*}|u|^{\alpha p^*}\, x\cdot DK(x)\biggr\}dx, 
\end{aligned}
\end{equation}
which yields to a contradiction when \eqref{non_ex} holds since the left hand side of \eqref{VI_final} is positive being $k\ge\alpha p^*$, $\lambda>0$ and thanks also to Lemma \ref{gprop}-{\bf d)}, while the right hand side of  \eqref{VI_final} is non positive.
\end{proof}

Moreover, in what follows we will prove the compactness of the inverse map of $G$ defined in \eqref{ggra} in some sense until $\alpha p^*$, taking into account similar results in \cite{DS, LWW3}.

\begin{theorem}\label{immc}
The map $|G^{-1}|^\alpha: D^{1,p}(\mathbb R^N)\to L^q(\mathbb R^N)$ is continuous for $p\le q \le p^*$ and
$|G^{-1}|^\alpha: D^{1,p}(\mathbb R^N)\to L^q_{loc}(\mathbb R^N)$ is compact for $p\le  q < p^*$. Moreover, it holds the following
\begin{equation}\label{weak_conv}
 G^{-1}(v_n(x))\to G^{-1}(v(x))\text { a.e. } x\in\mathbb{R}^N.
\end{equation}
\end{theorem}

\begin{proof}
Take $v\in D^{1,p}(\mathbb R^N)$ and, by using \eqref{alphap} and
 Sobolev's inequality, we have
\begin{equation}\label{dva}\|G^{-1}(v)\|_{\alpha p^*}=\||G^{-1}(v)|^\alpha\|_{p^*}^{1/\alpha}\le S^{-1/\alpha p}\|D(|G^{-1}(v)|^\alpha)\|_p^{1/\alpha}\le (2/S)^{1/\alpha p} \|Dv\|_p^{1/\alpha},
\end{equation}
where $S$ is Sobolev's constant, i.e.
$S=\inf\left\{\|D u\|^p_{p}\cdot\|u\|^{-p}_{p^{*}}: u\in D^{1,p}(\mathbb{R}^N)\setminus\{0\}\right\}.$
So that \eqref{dva} imply 
\begin{equation}\label{dva2}\||G^{-1}(v)|^\alpha\|_{p^*}=\|G^{-1}(v)\|_{\alpha p^*}^\alpha\le C \|Dv\|_p,
\end{equation}
so that we get continuity for $q= p^*$. To prove continuity for $q<p^*$, consider $x,y>1$ with $1<x<q$ and use  H\"older's inequality with exponents $y$ and $y'>1$
$$\begin{aligned}
 \|G^{-1}(v)\|_{\alpha q}=&\biggl(\int_{\mathbb R^N} \bigl(|G^{-1}(v)|^\alpha \bigr)^{q-x} \bigl(|G^{-1}(v)|^\alpha\bigr )^{x} dx\biggr)^{1/\alpha q}\\
&\le \biggl( \int_{\mathbb R^N}\bigl ( |G^{-1}(v)|^{\alpha}\bigr)^{(q-x)y'} dx \biggr)^{1/\alpha q y'}  \biggl(\int_{\mathbb R^N} (|G^{-1}(v)|^\alpha)^{xy} dx\biggr)^{1/\alpha q y}\end{aligned}$$
Now choose $y=p^*/x>1$, $y'=p^*/(p^*-x)$
so that we get
$$\begin{aligned}
 \|G^{-1}(v)\|_{\alpha q}&\le 
 \biggl( \int_{\mathbb R^N}\bigl ( |G^{-1}(v)|^{\alpha}\bigr)^{(q-x)p^*/(p^*-x)} dx \biggr)^{(p^*-x)/\alpha q p^*}  \biggl(\int_{\mathbb R^N} (|G^{-1}(v)|^\alpha)^{p^*} dx\biggr)^{x/\alpha q p^*}
\\&\le\|v\|_{(q-x)p^*/(p^*-x)}^{(q-x)/\alpha q}
\cdot \||G^{-1}(v)|^{\alpha}\|_{p^*}^{x/\alpha q }\le C \|Dv\|_{p}^{(q-x)/\alpha q} \cdot\|Dv\|_p^{x/\alpha q }=C \|Dv\|_{p}^{1/\alpha}, \end{aligned}$$
where we have used Lemma \ref{gprop}-{\bf g)},  \eqref{dva2} and the continuity of embedding $D^{1,p}(\mathbb R^N)\hookrightarrow L^{(q-x)p^*/(p^*-x)}(\mathbb R^N)$,  since $q<p^*$.
In turn $\||G^{-1}(v)|^\alpha\|_q\le C \|Dv\|_{p}$.

To prove the compactness of the map 
$|G^{-1}|^\alpha:D^{1,p}(\mathbb R^N)\hookrightarrow L^q_{loc}(\mathbb R^N)$ for $p\le  q < p^*$, 
we start from a bounded sequence $(v_n)_n$ in $D^{1,p}(\mathbb R^N)$, so that up to subsequences $v_n \rightharpoonup v$
in $D^{1,p}(\mathbb R^N)$
 then, by Lemma \ref{g-1N} we have $|G^{-1}(v_n)|^\alpha \rightharpoonup |G^{-1}(v)|^\alpha$ in $D^{1,p}(\mathbb R^N)$
and by the compactness of the embedding of $D^{1,p}(\mathbb R^N)$ in $L^s_{loc}(\mathbb R^N)$ for any $1<s<p^*$ we have
\begin{equation}\label{cbb}
|G^{-1}(v_n)|^\alpha\to |G^{-1}(v)|^\alpha \,\ \text{in} \,\ L^{s}(\omega),\quad \omega\Subset\mathbb{R}^N,\quad 1\le s< p^*.
\end{equation}
Finally, from \eqref{cbb}, by using an increasing sequence of compact sets whose union is $\mathbb R^N$
 and a diagonal argument, we get \eqref{weak_conv}.

\end{proof}

In order to prove the regularity of $F_\lambda$, we need to analyze the regularity of 
$$J(v)=\int_{\mathbb{R}^N} V|G^{-1}(v)|^{k}dx\quad \text{and} \quad 
H(v)=\int_{\mathbb{R}^N} K|G^{-1}(v)|^{\alpha p^{*}}dx$$

\begin{lemma}\label{lem1}
If $V\in L^{r}\left(\mathbb{R}^N\right)$ and $1<\alpha<k<\alpha p^*$, then $J(v)$
is weakly continuous on $D^{1,p}(\mathbb R^N)$.
 Moreover, $J$ is continuously differentiable and
$J': D^{1,p}(\mathbb R^N)\to [D^{1,p}(\mathbb R^N)]'$, for all $\psi\in D^{1,p}(\mathbb R^N)$, is given by
\begin{equation}\label{J-G}
J'(v)\psi=k\int_{\mathbb{R}^N} V\frac{|G^{-1}(v)|^{k-2}G^{-1}(v)}{g(G^{-1}(v))}\psi dx.
\end{equation}
\end{lemma}

\begin{proof}
For any $v\in D^{1,p}(\mathbb R^N)$, by Theorem \ref{immc}, then $|G^{-1}(v)|^\alpha\in  L^{p^{*}}(\mathbb{R}^N)$, so that by H\"older inequality with exponents $r=\alpha p^{*}/(\alpha p^{*}-k)$ and
$r'=\alpha p^{*}/k$ we have 
\begin{equation}\label{hol}
\|V|G^{-1}(v)|^k\|_1\le \|V\|_{r}\|G^{-1}(v)\|_{\alpha p^{*}}^{k}
\end{equation}
This implies that $J$ is well defined.
Let $(v_{n})_{n}\in D^{1,p}(\mathbb R^N)$ such that $v_n \rightharpoonup v$ in $D^{1,p}(\mathbb R^N)$,
thus, $(v_{n})_{n}$ is bounded in $D^{1,p}(\mathbb R^N)$ and, by Theorem \ref{immc}, also
 $(G^{-1}(v_n))_{n}$ is bounded in $L^{\alpha p^{*}}(\mathbb{R}^N)$ and $(|G^{-1}(v_n)|^k)_n$ in
$L^{\alpha p^{*}/k}(\mathbb{R}^N)$.

Since $V\in L^r(\mathbb R^N)$, by \eqref{weak_conv} and \eqref{hol}, the latter applied with $v_n$ instead of $v$, we get the weakly continuity of $J$, thanks to Lebesgue convergence Theorem, that is
\begin{equation}\label{wc}
J(v_n)=\int_{\mathbb{R}^N} V|G^{-1}(v_n)|^{k}dx\to \int_{\mathbb{R}^N} V|G^{-1}(v)|^{k}dx=J(v).
\end{equation}

In order to prove the Fr\'echet differentiability, that is $J\in C^{1}$, it is enough to show that $J$ is G\^ateaux differentiable and has a continuous G\^ateaux derivative on 
$D^{1,p}(\mathbb R^N)$.
First, consider $v, \psi\in D^{1,p}(\mathbb{R}^N)$ and $0<|t|<1$, so that
\begin{equation}\label{jvf}
\frac{J(v+t\psi)-J(v)}{t}=\int_{\mathbb R^N} V \frac{|G^{-1}(v+t\psi)|^{k}-|G^{-1}(v)|^{k}}{t}\, dx.
\end{equation}
Using the mean value theorem, there exists
$\delta\in(0,1)$ such that
\begin{equation}\label{gg}\begin{aligned}
\frac{\left| |G^{-1}(v+t\psi)|^{k}-|G^{-1}(v)|^{k}\right|}{|t|}&=k|G^{-1}(v+t\delta\psi)|^{k-1}|(G^{-1}(v+t\delta\psi))'||\psi|\\
&\le c|v+t\delta\psi|^{(k-\alpha)/\alpha}|\psi|\le c \left(|v|^{(k-\alpha)/\alpha}|\psi|+|\psi|^{k/\alpha}\right),
\end{aligned}\end{equation}
with $c>0$, where we have used, for the first inequality, the following condition
\begin{equation}\label{cond1}
 \frac{|G^{-1}(s)|^{k-1}(G^{-1}(s))'}{s^{(k-\alpha)/\alpha}}\le c
\end{equation}
which holds for on bounded sets in $\mathbb R_0^+$ by using Lemma \ref{gprop}-{\bf g), i)} .
While, the last inequality in \eqref{gg} follows from the elementary formula $(a+b)^{r}\le C( a^{r}+b^{r})$, $a,b,r>0$ and $C>0$.

Now, by applying H\"older's inequality twice with exponents $r$, $\alpha p^{*}/(k-\alpha)$, $p^{*}$ and $r$, $\alpha p^*/k$, we get
$$\int_{\mathbb{R}^N}V\left(|v|^{(k-\alpha)/\alpha}|\psi|+|\psi|^{k/\alpha}\right)dx
\le\left\|V\right\|_{r}\left\|\psi\right\|_{ p^{*}}\left(\left\|v\right\|_{ p^{*}}^{(k-\alpha)/\alpha}+\left\|\psi\right\|_{p^{*}}^{k/\alpha-1}\right).
$$
The right hand side  of the above inequality is finite thanks to the suitable summabilities of the functions $V,\psi,v$. Thus, by letting $t\to0$ in \eqref{jvf} using \eqref{gg}, from the
Lebesgue dominated convergence theorem, $J$ is G\^ateaux differentiable.

Before proving the continuity of the G\^ateaux derivative, we claim that the functional $J'(v)$ is well defined for every $v\in D^{1,p}(\mathbb{R}^N)$, that is  $J'(v)\in [D^{1,p}(\mathbb{R}^N)]'$. Indeed,
as in Lemma 2.5 in \cite{DS}, consider $\psi_n\to 0$ in $D^{1,p}(\mathbb R^N)$ so that, by \eqref{wc}, we get 
\begin{equation}\label{wc2}
\int_{\mathbb{R}^N}V|G^{-1}(\psi_n)|^{k} dx\to 0, \quad \text{as}\,\, n\to \infty
\end{equation}
Now, by using that $g\ge 1$, \eqref{G1a} and Lemma \ref{gprop}-{\bf i)} we have 
$$\begin{aligned}
&\int_{\mathbb{R}^N}V|G^{-1}(v)|^{k-1}G^{-1}(v)'\psi_n dx \\
&\le C_1 \biggl(\int_{\mathbb{R}^N}V|G^{-1}(v)|^{k} dx\biggr)^{1-1/k}\cdot \biggl(\int_{\mathbb{R}^N}V|G^{-1}(\psi_n)|^k dx\biggr)^{1/k}
\\&\quad+C_2\biggl(\int_{\mathbb{R}^N}V|G^{-1}(v)|^{k}\biggr)^{1-\alpha/k}\cdot
\biggl(\int_{\mathbb{R}^N}V|G^{-1}(\psi_n)|^k dx\biggr)^{\alpha/k}
\end{aligned}$$
where we used H\"older's inequality and the fact that $\alpha<k$. Now, applying \eqref{wc2}, we get 
$$\int_{\mathbb{R}^N}V|G^{-1}(v)|^{k-1}G^{-1}(v)'\psi_n dx\to0\quad \text{as}\,\, n\to \infty$$
obtaining the continuity of $J'$ and, consequently, the claim. 

Now, to get (Fr\'echet) differentiability  we 
check that the G\^ateaux derivative $J':D^{1,p}(\mathbb{R}^N)\to [D^{1,p}(\mathbb{R}^N)]'$ defined in \eqref{J-G} is continuous. To reach the claim,
consider $v_n\to v$ in $D^{1,p}(\mathbb R^N)$, 
so that there exists $U\in L^{p^{*}}(\mathbb{R}^N)$ such that
$|v_n|\le U(x)$ a.e. in $\mathbb R^N$. Now, define $W(v)=V|G^{-1}(v)|^{k-2}G^{-1}(v)(G^{-1}(v))'$ hence, by \eqref{cond1} and Young's inequality with exponents $\alpha(p^*-1)/(\alpha p^*-k)$ and $\alpha(p^*-1)/(k-\alpha)$, we get for $c>0$
$$\begin{aligned}
|W(v_n)-W(v)|^{(p^{*})'}&\le c ( |W(v_n)|^{(p^{*})'}+|W(v)|^{(p^{*})'})\\
&\le c\bigl(|V|^{p^*/(p^*-1)}|v_n|^{p^*(k-\alpha)/\alpha(p^*-1)}+|V|^{p^*/(p^*-1)}|v|^{p^*(k-\alpha)/\alpha(p^*-1)}\bigr)\\
&\le c\bigl(|V|^r +|v_n|^{p^*}
+|v|^{p^*}\bigr)\le c\bigl(|V|^r +U^{p^*}
+|v|^{p^*}\bigr)\in L^1(\mathbb R^N).
\end{aligned}$$
In turn, 
Lebesgue dominated convergence Theorem gives $\|W(v_n)-W(v)\|_{(p^{*})'}\to0$ as $n\to\infty$. Consequently,
by H\"older's inequality, for $\psi\in D^{1,p}(\mathbb R^N)$, we have
$$|\bigl(J'(v_n)-J'(v)\bigr)\psi|\le k\int_{\mathbb{R}^N} |W(v_n)-W(v)||\psi|dx\le k\|W(v_n)-W(v)\|_{(p^{*})'}
\|\psi\|_{p^*}\to 0,$$
as $n\to \infty$. Namely, $J\in C^{1}$.
\end{proof}


\begin{lemma}\label{lem2}
If $K\in L^{\infty}(\mathbb{R}^N)\cap C(\mathbb{R}^N)$, then $H$ is continuously differentiable in $D^{1,p}(\mathbb R^N)$ and
 its derivative $H': D^{1,p}(\mathbb R^N)\to [D^{1,p}(\mathbb R^N)]'$, for all $v, \psi\in D^{1,p}(\mathbb R^N)$, is given by
$$H'(v)\psi=\alpha p^{*}\int_{\mathbb{R}^N} K\frac{|G^{-1}(v)|^{\alpha p^{*}-2}G^{-1}(v)}{g(G^{-1}(v))}\psi dx,$$
.
\end{lemma}

\begin{proof}
The proof relies on the one of Lemma \ref{lem1} but with some adjustments. First, note that here there are no conditions on the exponent $k$.
Trivially, the functional is well defined and weakly continuous on $D^{1,p}(\mathbb R^N)$.
In order to prove the G\^ateaux differentiability of $H$, we use again the Mean Value Theorem but, instead of \eqref{gg}, we have
$$\begin{aligned}&\frac{\left| |G^{-1}(v+t\psi)|^{\alpha p^*}-|G^{-1}(v)|^{\alpha p^*}\right|}{|t|}=\alpha p^*|G^{-1}(v+t\delta\psi)|^{\alpha p^*-1}|(G^{-1}(v+t\delta\psi))'||\psi|\\
&\quad \le c |G^{-1}(v+t\delta\psi)|^{\alpha(p^*-1)}|\psi|\le c|v+t\delta\psi|^{p^*-1}|\psi|\le c \left(|v|^{p^*-1}|\psi|+|\psi|^{p^*}\right),
\end{aligned}$$
where we have used Lemma \ref{gprop}-{\bf g), i)} and the elementary formula $(a+b)^{r}\le C( a^{r}+b^{r})$, $a,b,r>0$ and $C>0$.
Concerning the well definition of $H$ on $[D^{1,p}(\mathbb R^N)]'$, we get 
$$\begin{aligned}
\int_{\mathbb{R}^N}&K|G^{-1}(v)|^{\alpha p^*-1}G^{-1}(v)'\psi_n dx \le C_1\int_{\mathbb{R}^N}K|G^{-1}(v)|^{\alpha p^*-1}G^{-1}(v)'G^{-1}(\psi_n) dx\\
&\quad+C_2\int_{\mathbb{R}^N}K|G^{-1}(v)|^{\alpha p^*-1}G^{-1}(v)'|G^{-1}(\psi_n)|^\alpha dx\\
&\le C_1 \biggl(\int_{\mathbb{R}^N}|G^{-1}(v)|^{\alpha p^*} dx\biggr)^{1-1/\alpha p^*}\cdot \biggl(\int_{\mathbb{R}^N}|G^{-1}(\psi_n)|^{\alpha p^*} dx\biggr)^{1/\alpha p^*}\\
&\quad+C_2\int_{\mathbb{R}^N}|G^{-1}(v)|^{\alpha( p^*-1)}|G^{-1}(\psi_n)|^\alpha dx\\
&\le C_1 \biggl(\int_{\mathbb{R}^N}|G^{-1}(v)|^{\alpha p^*} dx\biggr)^{1-1/\alpha p^*}\cdot \biggl(\int_{\mathbb{R}^N}|G^{-1}(\psi_n)|^{\alpha p^*} dx\biggr)^{1/\alpha p^*}\\
&\quad+C_2\biggl(\int_{\mathbb{R}^N}|G^{-1}(v)|^{\alpha p^*}\biggr)^{1-1/p^*}\cdot
\biggl(\int_{\mathbb{R}^N}|G^{-1}(\psi_n)|^{\alpha p^*} dx\biggr)^{1/p^*}
\end{aligned}$$
where, for the integral near $C_1$ we use that $g\ge 1$ and H\"older's inequality with exponents $\alpha p^*$ and $\alpha p^*/(\alpha p^*-1)$, while for the integral with $C_2$ we use  in Lemma \ref{gprop}-{\bf i)} and H\"older's inequality with exponents $p^*$ and $p^*/(p^*-1)$.

Finally, the continuity of the G\^ateaux derivative follows by defining $\tilde W$ instead of $W$ as $\tilde W(v)=|G^{-1}(v)|^{\alpha p^*-2}G^{-1}(v)(G^{-1}(v))'$ and  from the following
$$|W(v)|^{(p^*)'}\le |G^{-1}(v)|^{\alpha p^*}\le |v|^{p^*}$$
where Lemma \ref{gprop}-{\bf g), i)} is used.
\end{proof}

Now,
if we consider $v_n\to v$ in $D^{1,p}(\mathbb R^N)$, then by Lemma \ref{g-1N}, also $|G^{-1}(v_n)|^\alpha\to |G^{-1}(v)|^\alpha$ in $D^{1,p}(\mathbb R^N)$. 
Since the first term of $F_{\lambda}$ is a norm with exponent $p>1$,  and thanks to Lemmas \ref{lem1}
and \ref{lem2} with $1<\alpha<k$, then we get
$F_{\lambda}\in C^{1}(D^{1,p}(\mathbb R^N))$, and 
$F'_\lambda: D^{1,p}(\mathbb R^N)\to (D^{1,p}(\mathbb R^N))'$, for all $v, \psi\in D^{1,p}(\mathbb R^N)$, is given by
\begin{equation}\label{E'}
\begin{aligned}
F'_{\lambda}(v)\psi&=\int_{\mathbb{R}^N}|D v|^{p-2}D v D\psi dx 
-\lambda\int_{\mathbb{R}^N} V\frac{|G^{-1}(v)|^{k-2}G^{-1}(v)}{g(G^{-1}(v))}\psi dx\\
&\qquad\qquad-\beta\int_{\mathbb{R}^N} K\frac{|G^{-1}(v)|^{\alpha p^{*}-2}G^{-1}(v)}{g(G^{-1}(v))}\psi dx.
\end{aligned}
\end{equation}

We say  that $v\in D^{1,p}(\mathbb R^N)$ is a (weak) {\it solution} of equation \eqref{Psc} if
$$F'_{\lambda}(v)\psi=0\quad\text{ for all }\psi\in D^{1,p}(\mathbb R^N),$$

Clearly, (weak) {\it solutions} of \eqref{Psc} are exactly critical points of the Euler–Lagrange functional $E_{\lambda}$, or equivalently $F_\lambda$, associated with \eqref{Psc}.
Moreover, every critical point of $F_\lambda$ correspond to a solution of the following equation
\begin{equation}\label{npsc}
-\Delta_p v=\lambda V\frac{|G^{-1}(v)|^{k-2}G^{-1}(v)}{g(G^{-1}(v))}+\beta K\frac{|G^{-1}(v)|^{\alpha p^{*}-2}G^{-1}(v)}{g(G^{-1}(v))}.
\end{equation}
As a consequence,  Theorems \ref{t11} and \ref{T_E>0}, whose statements are given in the Introduction in terms of $u$ satisfying \eqref{Psc}, can be stated in terms of $v=G(u)$ solutions of \eqref{npsc}.

A key role in the proof of our results is the concentration compactness principles by Lions. For a detailed discussion on them, we refer to \cite{BBF} and \cite{BBFS}.
In particular, we are interested in the second concentration compactness principle, which regards a possible concentration only at finite points and
where two different types (since we are in unbounded domains) of convergences are considered: the  tight convergence
of measures, whose symbol is $\overset{\ast}{\rightharpoonup}$, and the "weak" convergence, denoted with  $\rightharpoonup$. Precisely,

\begin{lemma}\label{lem4.1} (Lemma \textnormal{I.1}, \textnormal{\cite{LionsRev1}})
Assume $\Omega\subset \mathbb R^N$ a domain, $1\le p< N$.
Let $(u_{n})_{n}$ be a bounded sequence in $D^{1,p}(\Omega)$ converging weakly to some $u$ and such that
$|Du_{n}|^{p}dx\rightharpoonup\mu$ and either $|u_{n}|^{p^*}dx\rightharpoonup\nu$ if $\Omega$ is bounded
or $|u_{n}|^{p^*}\overset{\ast}{\rightharpoonup}\nu$ if $\Omega$ is unbounded, where
$\mu$, $\nu$ are bounded nonnegative measures on $\Omega$. Then there exists some at most countable set
$J$ such that
$$\nu=|u|^{p^*}+ \sum_{j\in J}\nu_{j}\delta_{x_{j}}, \quad \mu\ge|Du|^{p}+\sum_{j\in J}\mu_{j}\delta_{x_{j}}, \qquad \nu_j\ge 0, \,\,\mu_j\ge 0 $$
with 
$\displaystyle{S\nu_{j}^{p/p^*}\le \mu_{j}}$ and 
$\displaystyle{\sum_{j\in J}\nu_{j}^{p/p^*}<\infty}$,
where $(x_{j})_{j\in J}$ are distinct points in $\Omega$, $\delta_{x}$ is the Dirac-mass of mass 1 concentrated at
$x\in\Omega$.
\end{lemma}

Since the tight convergence excludes a possible  concentration at infinity, by using the lemma above, in order to get compactness, it remains only to show that concentration around points, described by $\nu_j$, cannot occur. 
However, the proof of tightness by definition as well as using the first concentration compactness principle leads to rather cumbersome and tricky calculations. 
To overcome these difficulties, Chabrowskii presented a version at infinity of the
second principle, cfr. Proposition 2  in \cite{cha95ca}, see also Bianchi et al. in \cite{BCS} for the Laplacian case.
In this principle Chabrowskii manages to enclose the concentration at infinity in the parameter $\nu_\infty$, in according to the $\nu_j$s in Lemma \ref{lem4.1} so that the non concentration at infinity occurs if one proves that $\nu_\infty=0$.
Later Ben-Naoum et al. in \cite{BTW} obtain a version for the $p$-Laplacian, here reported for completeness.

\begin{proposition}\label{lemben} (Proposition \textnormal{3.3}, \textnormal{\cite{BTW}})
Let $(u_n)_n$ be a bounded sequence in $D^{1,p}(\mathbb R^N)$ and define
$$\nu_\infty=\lim_{R\to\infty} \limsup_{n\to\infty} \int_{|x|>R} |u_n|^{p^*} dx,\qquad
\mu_\infty=\lim_{R\to\infty} \limsup_{n\to\infty} \int_{|x|>R} |Du_n|^{p} dx.$$
Then, the quantities $\nu_\infty$ and $\mu_\infty$ exist and satisfy
$$\limsup_{n\to\infty} \int_{\mathbb R^N} |u_n|^{p^*} dx=\int_{\mathbb R^N} d\nu +\nu_\infty,\qquad \limsup_{n\to\infty} \int_{\mathbb R^N} |Du_n|^{p} dx=\int_{\mathbb R^N} d\mu +\mu_\infty,$$
with $S\nu_{\infty}^{p/p^*}\le \mu_\infty,$
where $\nu$ and $\mu$ are as in $(i)$ and $(ii)$ in Lemma \ref{lem4.1} and such that $(iii)$ holds.
\end{proposition}

Another key tool in proving the multiplicity result of solutions with positive energy, namely Theorem \ref{T_E>0}, is the Fountain Theorem. However, before its statement, we briefly report the setting needed. For some well-known basic definitions of actions, invariant functions, we refer to \cite{BBFS}.

\begin{enumerate}
\item[{\bf (A1)}] The compact group $\mathcal G$ acts isometrically on the space
$M=\overline{\bigoplus_{j\in\mathbb{N}}M_{j}}$, which is
a Banach space, where the spaces $M_{j}$ are $\mathcal G$-invariant and there exists a finite dimensional space $W$ such that, for
every $j\in \mathbb{N}$, $M_{j}\simeq W$ and the action of $\mathcal G$ on $W$ is admissible.
\end{enumerate}

From the decomposition of the Banach space $M$ in {\bf(A1)}, we define $Y_m$ and $Z_m$ as follows
\begin{equation}\label{mare}
Y_{m}:=\bigoplus_{j=0}^{m}M_{j}, \,\,\ Z_{m}:=\overline{\bigoplus_{j=m}^{\infty}M_{j}}
\end{equation}
and set $B_{m}:=\left\{u\in Y_{m}: \| u\| \le \rho_{m}\right\}$, $N_{m}=\left\{u\in Z_{m}: \|u\|=r_{m}\right\}$ where $\rho_{m}>r_{m}>0$.

Now we are ready to state the Fountain Theorem.

\begin{theorem}\textnormal{(Theorem 3.6, \textnormal{\cite{W}})}\label{FT}
Under assumption \textnormal{\bf(A1)}. Let $\varphi\in C^{1}(M,\mathbb{R})$ be an invariant functional.
If, for every $m\in\mathbb{N}$, there exists $\rho_{m}>r_{m}>0$ such that
\begin{enumerate}
\item[{\bf (A2)}] $a_{m}=\max_{u\in Y_m, \|u\|=\rho_m} \varphi(u)\le 0$,
\item[{\bf (A3)}] $b_{m}=\inf_{u\in Z_m, \|u\|=r_m} \varphi(u)\to\infty, \,\ m\to\infty$,
\item[{\bf (A4)}] $\varphi$ satisfies the (PS)$_{c}$ condition for every $c>0$,
\end{enumerate}
where $Y_{m}$ and $Z_{m}$ as in \eqref{mare}.
 Then $\varphi$ has an unbounded sequence of critical values.
\end{theorem}

\begin{remark}\label{aft}
In our setting, as done in \cite{BBFS}, we set $\mathcal G=\mathbb{Z}/2$, $M=D^{1,p}_T(\mathbb R^N)$ so that, since $D^{1,p}_T(\mathbb R^N)$ is a separable Banach space, there is a linearly independent sequence $(e_{j})_j$ such that the decomposition
in {\bf(A1)} holds with $M_{j}=X_j:=span\left\{e_{j}\right\}$. Note that $X_j$ are trivially $\mathcal{G}$-invariant and
isomorphic to $\mathbb R$. Thus, condition {\bf(A1)} is satisfied with $W=\mathbb R$.
\end{remark}

Finally, a crucial ingredient when a symmetric setting is involved is the following principle of symmetric criticality due to Palais, cfr. \cite{Palais}, \cite{PRRbook}, which states that  any critical point of $F_\lambda$ restricted on $D^{1,p}_T(\mathbb R^N)$, is a  critical point of the same functional on  $D^{1,p}(\mathbb R^N)$.
For further details, we refer to \cite{BBFS}.

\begin{lemma}
Let $v\in D^{1,p}_T(\mathbb R^N)$. If $F_\lambda(v)\psi=0$ for all $\psi\in D^{1,p}_T(\mathbb R^N)$, then $F_\lambda(v)\psi=0$ for all $\psi\in D^{1,p}(\mathbb R^N)$.
\end{lemma}

\section{On Palais Smale sequences}\label{psss}

We start this section by giving the definition of (PS)$_c$ sequence. 

\begin{definition}\label{PS}
Let $Y$ be a Banach space and $E: Y\to\mathbb{R}$ be a differentiable functional. A sequence
$(u_{n})_{n}\subset Y$ is called a (PS)$_{c}$ sequence for $E$ if $E(u_{k})\to c$ and $E'(u_{k})\to 0$
as $k\to\infty$.
Moreover, we say that $E$ satisfies the (PS)$_{c}$ condition if every (PS)$_{c}$ sequence for $E$ has a converging subsequence in $Y$.
\end{definition}

As standard, we need to deal with bounded (PS)$_{c}$ sequences for the functional $F_{\lambda}$ in \eqref{F}. In addition, a useful inequality holds only if $\alpha<k<\alpha p$ is proved in the next lemma. 

\begin{lemma}\label{lem5}
Assume $\alpha<k<\alpha p^*$ and the further condition $K\ge0$ in $\mathbb R^N$ if $p<k<\alpha p^*$. Let \eqref{V(x)1} and \eqref{K(x)1} be verified and consider $(v_n)_{n}\subset D^{1,p}(\mathbb R^N)$  a
(PS)$_{c}$ sequence for
$F_{\lambda}$ for all $c\in\mathbb{R}$.
Then $(v_n)_{n}$ is bounded in $D^{1,p}(\mathbb R^N)$.

In particular, if $\alpha<k<\alpha p$ and $c<0$, it holds
\begin{equation}\label{bound_u_np*}
\|D v_n\|_p <C_* \lambda^{\alpha/(\alpha p-k)}, \qquad C_*=\biggl[\frac{p(\alpha p^*-k)}{\alpha k(p^*-p)}\cdot\biggl(\frac{2}{S}\biggl)^{{k}/{\alpha p}} \|V\|_r\biggr]^{\alpha/(\alpha p-k)},
\end{equation}
where $S$ is Sobolev's constant.
\end{lemma}

\begin{proof}
We follow Lemma 4 in \cite{BBF}. Let $(v_{n})_{n}\subset D^{1,p}(\mathbb R^N)$ be a (PS)$_{c}$ sequence of $F_{\lambda}$ for all $c\in\mathbb{R}$
 that is, using Definition \ref{PS},
$F_{\lambda}(v_{n})=c+o(1)$, $F'_{\lambda}(v_{n})\psi=o(1)\|\psi\|$ as $n\to\infty$,
for every $\psi\in D^{1,p}(\mathbb R^N)$. Now take $\psi=G^{-1}(v_n)g(G^{-1}(v_n))$ as a test function, since $\psi\in D^{1,p}(\mathbb R^N)$ thanks to Remark \ref{gd1p}, and using $|D[G^{-1}(v_n)g(G^{-1}(v_n))]|\le\alpha |Dv_n|$ by \eqref{eq1}, we have
\begin{equation}\label{eprimo}\begin{aligned}
o(1)\|\psi\|&=F'_\lambda(v_n)(\psi)\\
&\le \alpha \|D v_n\|_p^p-\lambda \int_{\mathbb R^N} V|G^{-1}(v_n)|^k dx-\beta\int_{\mathbb R^N} K |G^{-1}(v_n)|^{\alpha p^*} dx
\end{aligned}\end{equation}
Now we disjoint the proof in two cases.

{\it Case $\alpha<k<\alpha p$:} using \eqref{eprimo}, thanks to Lemma \ref{gprop}-{\bf c), e)}, \eqref{dva} and
H\"older's inequalities with exponents $r$ and $r'$ we get
\begin{equation}\label{E-E'u}\begin{aligned}
c &+o(1) +o(1)\|Dv_{n}\|_p\ge F_{\lambda}(v_n)- \frac{1}{\alpha p^{*}}F'_\lambda(v_n)(G^{-1}(v_n)g(G^{-1}(v_n)))\\
&\ge \left(\frac{1}{p}-\frac{1}{ p^{*}}\right)\|D v_n\|_p^p-\lambda\left(\frac1k-\frac{1}{\alpha p^{*}}\right)\|V\|_r\|G^{-1}(v_n)\|^{k}_{\alpha p^*}\\
&\ge \left(\frac{1}{p}-\frac{1}{ p^{*}}\right)\|D v_n\|_p^p-\lambda\left(\frac1k-\frac{1}{\alpha p^{*}}\right)(2/S)^{k/\alpha p}\|V\|_r\|D v_n\|^{k/\alpha}_{p}.
\end{aligned}\end{equation}
Thus, since $k<\alpha p <\alpha p^*$, we conclude that
$\|Dv_n\|_p$ should be bounded.

{\it Case $\alpha p< k<\alpha p^*$:} arguing  as in \eqref{E-E'u}, with $1/\alpha p^*$ replaced by $1/k$,
since $K(x)\geq 0$ in $\mathbb R^N$, we obtain
$$
\begin{aligned}
c &+o(1) +o(1)\|Dv_{n}\|_p\ge F_{\lambda}(v_n)- \frac{1}{k}F'_\lambda(v_n)(G^{-1}(v_n)g(G^{-1}(v_n)))\\
&\ge \left(\frac{1}{p}-\frac{\alpha}{k}\right)\|D v_n\|_p^p-\beta\left(\frac{1}{\alpha p^*}-\frac{1}{k}\right)
\int_{\mathbb R^N} K |G^{-1}(v_n)|^{\alpha p^*} dx \ge \left(\frac{1}{p}-\frac{\alpha}{k}\right)\|D v_n\|_p^p.\end{aligned}$$
The conclusion follows immediately since $k>\alpha  p$. 
For the proof of inequality \eqref{bound_u_np*}, as in Lemma 4 in \cite{BBF}, it is enough to observe that if $c<0$, then \eqref{E-E'u} gives
$$\|D v_n\|_p^{(\alpha p -k)/\alpha} < \,\frac{\lambda p(\alpha p^* -k)}{\alpha k (p^*-p)}\|V\|_r (2/S)^{k/\alpha p},
$$
from which \eqref{bound_u_np*} follows immediately.
\end{proof}

In the lemma below we prove the validity of (PS)$_{c}$ condition of $F_\lambda$, that is the point in which the lack of compactness becomes manifest.

\begin{lemma}\label{lem6}
Suppose \eqref{V(x)1}, \eqref{K(x)1}, $\alpha<k<\alpha p$, $c<0$ and define $c_1$, $c_2$ as follows
\begin{equation}\label{c12}
c_1= \frac{1}{N}\biggl(\frac{S}{2}\biggr)^{p^*/(p^*-p)}, \qquad c_2=\biggl[\frac{2p(\alpha p^*-k)}{S\alpha (p^*-p)}\biggr]^{k/(\alpha p-k)} \frac{1}{k^{\alpha p/(\alpha p-k)}}.
\end{equation}
Then 
\begin{enumerate}
    \item[(I)] For any $\lambda>0$, there exists $\beta^*_{PS}>0$ defined as follows
    \begin{equation}\label{betastar_up}
\beta^*_{PS}=\biggl(\frac{c_1}{c_2}\biggr)^{(p^*-p)/p}\cdot\frac{1}{ (\lambda \|V\|_r)^{\alpha (p^*-p)/(\alpha p-k)}}\cdot\frac{1}{\|K\|_\infty},
\end{equation}
    such that for every $\beta\in(0,\beta^*_{PS})$, then $F_{\lambda}$ satisfies (PS)$_{c}$ condition.
    \item[(II)] For any $\beta>0$, there exists $\lambda^*_{PS}>0$ defined as follows
   \begin{equation}\label{lambdastar_up}
\lambda^*_{PS}=\biggl(\frac{c_1}{c_2}\biggr)^{(\alpha p-k)/\alpha p} \cdot \frac{1}{(\beta \|K\|_\infty)^{(\alpha p-k)/\alpha(p^*-p)}}\cdot\frac{1}{\|V\|_r},
\end{equation}
    such that for every $\lambda\in(0,\lambda^*_{PS})$, then $F_{\lambda}$ satisfies (PS)$_{c}$ condition.
\end{enumerate}
\end{lemma}

\begin{proof}
Let $(v_{n})_n$ be a (PS)$_{c}$ sequence, by Lemma \ref{lem5}, then $(v_{n})_n$ is bounded in $D^{1,p}(\mathbb R^N)$ and by Banach-Alaoglu's Theorem, there exists $v\in D^{1,p}(\mathbb R^N)$ such that, up to subsequences, we get $v_{n}\rightharpoonup v$ in $D^{1,p}(\mathbb R^N)$.
On the other hand, By Lemma \ref{g-1N}, follows that
$|G^{-1}(v_n)|^\alpha\rightharpoonup |G^{-1}(v)|^\alpha$ in $D^{1,p}(\mathbb R^N)$, so that $(|G^{-1}(v_n)|^\alpha)_n$ is bounded in $D^{1,p}(\mathbb R^N)$. In addition, Theorem \ref{immc} gives \eqref{weak_conv}. 
Applying in Proposition \ref{lemben}, there exist $\mu, \nu$, $\nu_\infty$, $\mu_\infty$ bounded nonnegative measures on $\mathbb R^N$ such that
\begin{equation}\label{nu}
\limsup_{n\to\infty} \int_{\mathbb R^N} |G^{-1}(v_n)|^{\alpha p^*} dx=\int_{\mathbb R^N} d\nu +\nu_\infty,
\end{equation}
and
$$\limsup_{n\to\infty} \int_{\mathbb R^N} \bigl|D(|G^{-1}(v_n)|^\alpha)\bigr|^{p} dx=\int_{\mathbb R^N} d\mu +\mu_\infty,$$
where
$$\nu_\infty=\lim_{R\to\infty} \limsup_{n\to\infty} \int_{|x|>R} |G^{-1}(v_n)|^{\alpha p^*} dx,\quad
\mu_\infty=\lim_{R\to\infty} \limsup_{n\to\infty} \int_{|x|>R} \bigl|D(|G^{-1}(v_n)|^\alpha)\bigr|^{p} dx.$$
Moreover, there exists at most countable set $J$, a family $(x_j)_{j\in J}$ of distinct points in
$\mathbb R^N$ and two families $(\nu_j)_{j\in J}, \,(\mu_j)_{j\in J}\in ]0,\infty[$ so that
$$\nu=|G^{-1}(v)|^{\alpha p^*}+ \sum_{j\in J}\nu_{j}\delta_{x_{j}}, \quad \nu_j\ge 0, \qquad\mu\ge\bigl|D(|G^{-1}(v_n)|^\alpha)\bigr|^{p}+\sum_{j\in J}\mu_{j}\delta_{x_{j}}, \quad \mu_j\ge 0,$$
satisfying
\begin{equation}\label{6.22}
S\nu_{j}^{p/p^{*}}\le\mu_{j}, \qquad S\nu_{\infty}^{p/p^*}\le \mu_\infty.
\end{equation}

Following an idea in \cite{BFsc}, where now we use \eqref{alphap} and that $V\in L^r(\mathbb R^N)$ with $r$ is given in \eqref{V(x)1}, it holds
\begin{equation}\label{6.23}
2\beta K(x_{j})\nu_{j}\ge \mu_{j}.
\end{equation}
Inequality \eqref{6.23} establishes that concentration of the measure $\mu$ cannot occur at points in which $K(x_{j})\le0$ being the right hand side positive. 
Consequently, $K(x_j)>0$ for all $j\in J$.
We claim that $J=\emptyset$. Indeed,
combining \eqref{6.22}$_1$ and \eqref{6.23}, we arrive to
\begin{equation}\label{J2}
 \nu_{j}\ge \left( \frac{ S}{2\beta K(x_{j})}\right)^{p^*/(p^*-p)}\ge \left( \frac{ S}{2\beta \|K\|_\infty}\right)^{p^*/(p^*-p)},\qquad j\in J.
\end{equation}
To reach the claim, we show that \eqref{J2} cannot occur for $\lambda$ or $\beta$ belonging to a suitable interval.
As in \cite{BBF}, assumption
\eqref{J2} forces that $|J|<\infty$.
Now, being $(v_{n})_n$ a (PS)$_c$ sequence,
choose again $G^{-1}(v_n)g(G^{-1}(v_n))\psi_{\varepsilon}$ as a test function in \eqref{E'}, where $\psi_{\varepsilon}(x)=
\psi\left((x-x_{j})/\varepsilon\right)$ for $0<\varepsilon<1$ and $\psi\in C_{c}^{\infty}(\mathbb{R}^N)$ such that $0\le\psi\le1$ in $\mathbb{R}^N$,
$\psi=0$ for $|x|>1$, $\psi=1$ for $|x|\le 1/2$. 
Then, using \eqref{eq1} and $0\le \psi_\varepsilon\le 1$, we have, for $n\to\infty$,
$$\begin{aligned}
c&+o(1)\|Dv_{n}\|_p= F_{\lambda}(v_{n})-\frac{1}{\alpha p^*} F'_{\lambda}(v_{n})(G^{-1}(v_n)g(G^{-1}(v_n)))\psi_{\varepsilon}\\
&= \frac{1}{p}\|Dv_{n}\|_p^p-\frac{\lambda}{k}
\int_{\mathbb{R}^N}V|G^{-1}(v_n)|^{k}dx-\frac1{\alpha p^*}\int_{\mathbb{R}^N} |Dv_n|^{p-1}\bigl| D(G^{-1}(v_n)g(G^{-1}(v_n)))\bigr|\psi_{\varepsilon} dx \\
&\qquad-\frac1{\alpha p^*}\int_{\mathbb{R}^N} G^{-1}(v_n)g(G^{-1}(v_n))|Dv_n|^{p-2}Dv_n \cdot D\psi_{\varepsilon} dx +\frac{\lambda}{\alpha p^*}\int_{\mathbb{R}^N}V|G^{-1}(v_{n})|^{k}\psi_{\varepsilon}dx\\
&\ge\frac{1}{N}\int_{\mathbb{R}^N}|Dv_n|^{p}\psi_\varepsilon dx -\frac{\lambda}{k}\int_{\mathbb R^N}V|G^{-1}(v_{n})|^k dx\\
&\qquad-\frac1{\alpha p^*}\int_{B_\epsilon(x_j)\setminus B_{\epsilon/2}(x_j)} G^{-1}(v_n)g(G^{-1}(v_n))|Dv_n|^{p-1}|D\psi_{\varepsilon}| dx.
\end{aligned}$$
In particular, since
$|D\psi_\varepsilon|\le  C/\varepsilon$ in the entire $\mathbb R^N$, by Lemma \ref{gprop}-{\bf e)}, Lemma \ref{lem5} and using H\"older's inequality twice with exponents $p$, $p'$ and $N/(N-p)$, $N/p$ respectively, we obtain
$$\begin{aligned}\biggl|\int_{B_\epsilon(x_j)\setminus B_{\epsilon/2}(x_j)}& G^{-1}(v_n)g(G^{-1}(v_n))|Dv_n|^{p-1}|D\psi_{\varepsilon}| dx\biggr|
\\& \le \alpha \int_{B_\epsilon(x_j)\setminus B_{\epsilon/2}(x_j)}
|v_n||Dv_n|^{p-1}|D\psi_{\varepsilon}| dx
\\&\le \alpha \biggl(\int_{\mathbb{R}^N}|Dv_n|^{p}\biggr)^{(p-1)/p}\cdot 
\biggl(\int_{B_\epsilon(x_j)\setminus B_{\epsilon/2}(x_j)}
|v_n|^p|D\psi_{\varepsilon}|^p dx\biggr)^{1/p}
\\&\le C 
\biggl(\int_{B_\epsilon(x_j)\setminus B_{\epsilon/2}(x_j)}
|v_n|^{p^*} dx\biggr)^{1/p^*}
\biggl(\int_{B_\epsilon(x_j)\setminus B_{\epsilon/2}(x_j)}
|D\psi_{\varepsilon}|^N dx\biggr)^{1/N}\\&\le C \biggl(\int_{B_\epsilon(x_j)\setminus B_{\epsilon/2}(x_j)}
|v_n|^{p^*} dx\biggr)^{1/p^*}\to 0,
\end{aligned}$$
as $\varepsilon\to0$, since $v_n\in L^{p^*}(\mathbb R^N).$
In turn, applying H\"older's inequality with exponents $r$, $r'$ and \eqref{dva} we arrive to 
\begin{equation}\label{ppstar}
c+o(1)\|Dv_{n}\|_p\ge\frac{1}{N}\int_{\mathbb{R}^N}|Dv_n|^{p} \psi_\varepsilon dx -\frac{\lambda}{k}\biggl(\frac{2}{S}\biggr)^{k/\alpha p}\|V\|_{r}\|Dv_{n}\|_{p}^{k/\alpha}+o(1)
\end{equation}
as $n\to\infty$.
Now, thanks to \eqref{alphap} and \eqref{bound_u_np*}, from \eqref{ppstar}, we get
$$c+o(1)\|Dv_{n}\|_p
> \frac{1}{2N}\int_{B_{\varepsilon/2}(x_j)}|D(|G^{-1}(v_{n})|^\alpha)|^{p} dx\,-\frac{(C_*)^{k/\alpha}}{k}\cdot \biggl(\frac{2}{S}\biggr)^{ k/\alpha p}\|V\|_{r}\lambda^{\alpha p/(\alpha p-k)}+o(1),
$$
where $C_*$ is given in \eqref{bound_u_np*}, so that, letting $n\to\infty$, $\varepsilon\to0$ and using \eqref{6.22} and \eqref{J2}, we arrive to
$$0>c>  c_1
\bigl(\beta\|K\|_{\infty}\bigr)^{-p/(p^*-p)}- c_2\bigl(\|V\|_r\lambda\bigr)^{\alpha p/(\alpha p-k)},$$
where we have replaced the value of $C_*$ and $c_1$, $c_2$ are defined in \eqref{c12}.
To obtain the required contradiction we need to have 
\begin{equation}\label{finee}
c_1>c_2 \bigl(\beta\|K\|_{\infty}\bigr)^{p/(p^*-p)}\bigl(\|V\|_r\lambda\bigr)^{\alpha p/(\alpha p-k)}.
\end{equation}
Consequently, since $k<\alpha p$, if we choose any $\beta>0$, then there exists 
$\lambda^*_{PS}$, defined in \eqref{lambdastar_up}, such that for every $\lambda\in(0,\lambda^*_{PS}]$, inequality \eqref{finee} is verified. Similarly, for any $\lambda>0$ fixed, there exists  $\beta^*_{PS}$, defined in \eqref{betastar_up}, such that for every $\beta\in(0,\beta^*_{PS}]$, inequality \eqref{finee} holds.
Thus $J=\emptyset$, concluding the proof of the claim.

On the other hand, following the idea of Chabrowski in \cite{cha95ca} and Ben-Naoum et. al in \cite{BTW}, also a possible concentration at infinity is refused.

Consequently, \eqref{nu} gives
$$\lim_{n\to\infty}\int_{\mathbb{R}^N}|G^{-1}(v_{n})|^{\alpha p^{*}}dx=\int_{\mathbb{R}^N}|G^{-1}(v)|^{\alpha p^{*}}dx.$$
Furthermore, since $G^{-1}(v_{n})(x)\to G^{-1}(v)(x)$ a.e. in $\mathbb R^N$ from \eqref{weak_conv}, then Brezis Lieb Lemma in
 \cite{BrezLieb}, implies
\begin{equation}\label{g-11}
\|G^{-1}(v_{n})\|^{\alpha p^{*}}_{\alpha p^{*}}-\|G^{-1}(v_{n})-G^{-1}(v)\|^{\alpha p^{*}}_{\alpha p^{*}}
= \|G^{-1}(v)\|^{\alpha p^{*}}_{\alpha p^{*}}+o(1).
\end{equation}
in other words, $\lim_{n\to\infty}\|G^{-1}(v_{n})-G^{-1}(v)\|_{\alpha p^{*}}=0$.
Using that by weak continuity of the functional $F_\lambda(v_n)\to F_\lambda(v)=c$, that is 
$$\begin{aligned}
\int_{\mathbb R^N}&|D v_n|^{p} dx-\int_{\mathbb R^N}|D v|^{p} dx +o(1)\\&=
\frac{\lambda p}{k}\int_{\mathbb{R}^N}V\bigl[|G^{-1}(v_n)|^{k}-|G^{-1}(v)|^{k}\bigr]dx
+\frac{\beta p}{\alpha p^{*}}\int_{\mathbb{R}^N}K\bigl[|G^{-1}(v_n)|^{\alpha p^{*}}
-|G^{-1}(v)|^{\alpha p^{*}}\bigr]dx
\end{aligned}$$
as $n\to\infty$. Thus, by \eqref{g-11}, the right hand side tends to 0, so that 
$\|Dv_n\|_p\to\|Dv\|_p$ as $n\to\infty$. Consequently, since $D^{1,p}(\mathbb R^N)$ is uniformly convex and by $Dv_n \rightharpoonup Du$ in $\mathbb R^N$, by Proposition 3.32 in \cite{Brezis}, we immediately get
$$\lim_{n\to\infty}\int_{\mathbb{R}^N}|D(v_n-v)|^{p}dx=0,$$
that is the strong convergence in $L^p(\mathbb R^N)$ of the sequence $(Dv_n)_n$.
\end{proof}

\section{The truncated functional $F_\infty$}\label{truncated}

In this section, we introduce $F_\infty$, the truncated functional  of $F_\lambda$, which has the very useful property to be bounded from below, differently from $F_\lambda$.

Taking into account Theorem \ref{immc}, in particular \eqref{dva2}, H\"older's inequality and by using \eqref{F}, for all $v\in D^{1,p}(\mathbb R^N)$ we have

$$F_{\lambda}(v)\ge \frac{1}{p}\| v\|^{p}-\lambda c_V\|v\|^{k/\alpha}-
\beta c_K\|v\|^{p^{*}}.$$
where $c_V=(2/S)^{k/\alpha p}\|V\|_r/k$ and $c_K=(2/S)^{p^*/p}\|K\|_\infty/\alpha p^*$ are positive constants.

\begin{figure}
     \centering
     \begin{subfigure}[b]{0.2\textwidth}
         \centering
         \includegraphics[width=\textwidth]{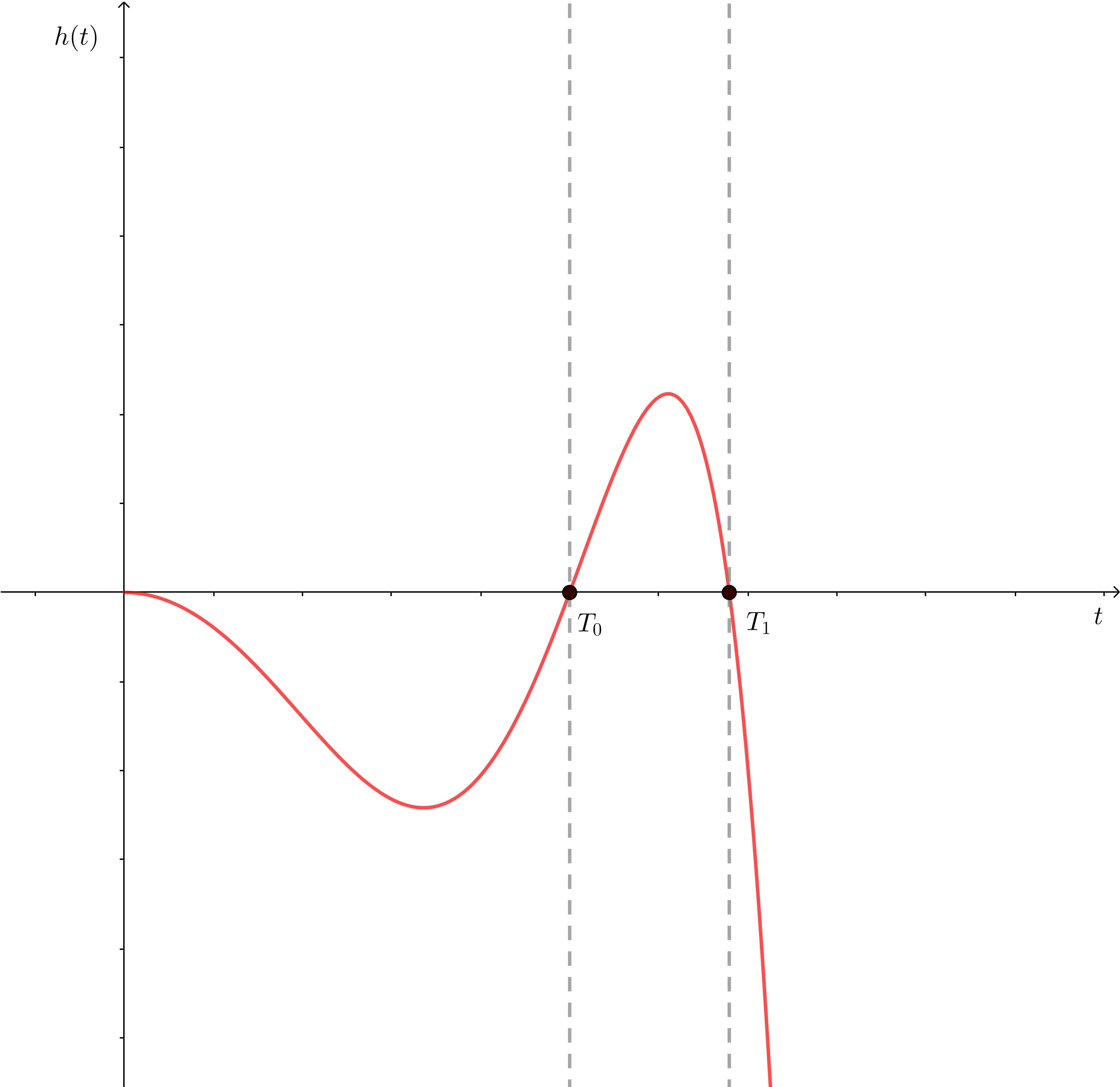}
         \caption{$ h$ if $\hat{h}(T_M)> 0$ }
         \label{h1}
     \end{subfigure}
     \hfill
     \begin{subfigure}[b]{0.2\textwidth}
         \centering
         \includegraphics[width=\textwidth]{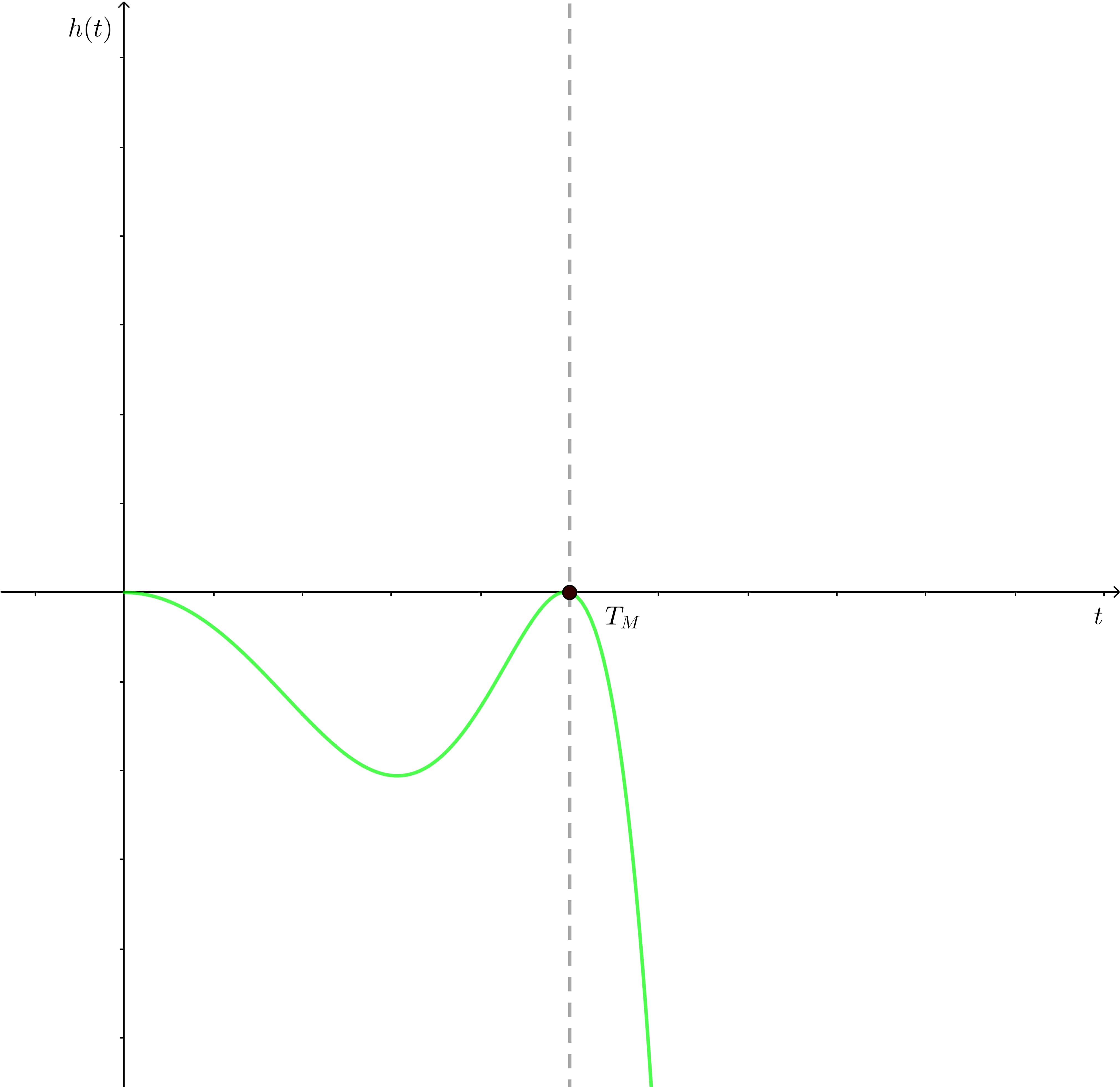}
         \caption{$ h$ if $\hat{h}(T_M)= 0$ }
         \label{h2}
     \end{subfigure}
     \hfill
     \begin{subfigure}[b]{0.2\textwidth}
         \centering
         \includegraphics[width=\textwidth]{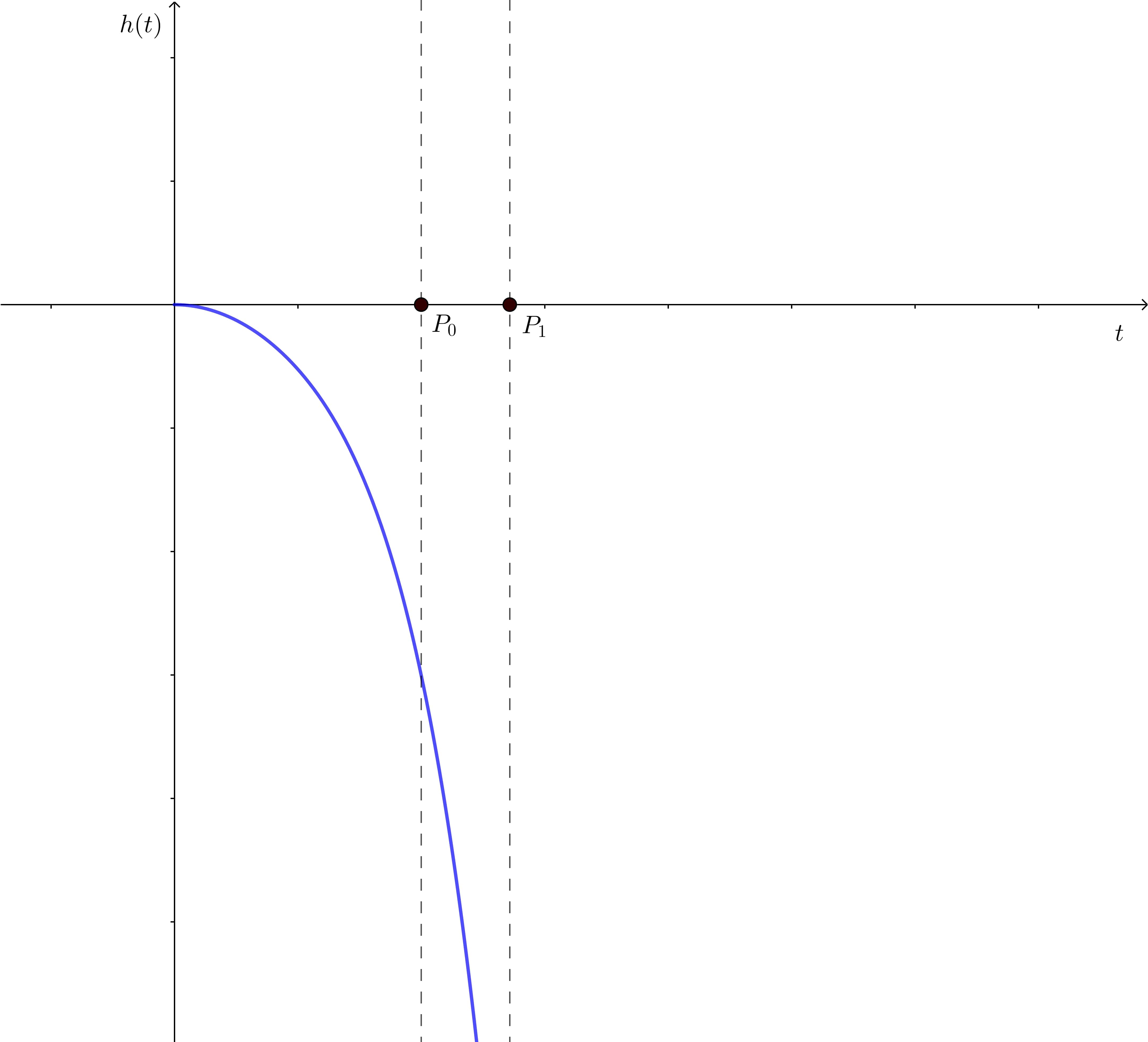}
         \caption{$ h$ if $\hat{h}(T_M)< 0$ }
         \label{h3}
     \end{subfigure}
        \caption{$h(t)$}
        \label{fig:three graphs}
\end{figure}

For simplicity, just to describe some qualitative properties of $F_\lambda$, we  define the auxiliary function $h(t)=t^{p}/p-\lambda c_Vt^{k/\alpha}-\beta c_Kt^{p^{*}}$ in $\mathbb R_0^+$. 
Indeed, even if $h$ is smaller than $F_\lambda$, it will appear in the proof of Theorem \ref{1.1} that the behaviour of $F_\lambda$ will be roughly the same of $h$ near $0$.
Following the same argument of Section 4 in \cite{BFsc}, we write 
$$h(t)=t^{ k/\alpha}\hat{h}(t),\qquad \hat{h}(t):= -\lambda c_V + \frac 1pt^{(\alpha p-k)/\alpha }-\beta c_Kt^{(\alpha p^{*}-k)/\alpha},$$
with $\hat{h}(0)<0$, $\hat{h}(t)\to-\infty$ as $t\to\infty$ and there is a unique point $T_M>0$ such that
$$\hat{h}'(T_M)=0, \quad T_M=\biggl[\frac{\alpha p-k}{\beta c_Kp(\alpha p^*-k)}\biggr]^{1/(p^*-p)}, \quad   \hat{h}'(t)>0 \text{ for } 0<t<T_M.$$
Hence, the maximum value of the function $\hat{h}$ is given by
\begin{equation}\label{hT}
\hat{h}(T_M)=\frac {\alpha(p^*-p)}{\alpha p-k}\,\bigl(\beta c_K\bigr)^{-(\alpha p-k)/\alpha(p^*-p)}\biggl(\frac{\alpha p-k}{p(\alpha p^*-k)}\biggr)^{(\alpha p^*-k)/\alpha(p^*-p)}-\lambda c_V.\end{equation}
In turn, since $\hat h$ and $h$ have the same zeros and the same sign in  $\mathbb R^+$,
if $\hat{h}(T_M)\ge 0$,  then there exists $T_0\in(0,T_M) $ such that $\hat h(T_0)=h(T_0)=0$ so that $\hat{h}(t)<0$ in  $(0,T_0)$ and consequently $h(t)<0$ in $(0,T_0)$, cfr. Figure \ref {h1} and Figure \ref{h2}.

While,  if $\hat{h}(T_M)<0$ then $\hat{h}(t),\, h(t)<0$ in  $\mathbb R^+$, so that  $h$ does not have zeros in all of $\mathbb R^+$. In addition, in this latter case $h$ is nonincreasing, indeed 
$$h'(t)=t^{k/\alpha}\hat h(t)\biggl[\frac k\alpha\frac 1t+\frac{\hat h'(t)}{\hat h(t)}\biggr]=
t^{k/\alpha}\hat h(t)\biggl(\log\mathfrak{Z}(t)\biggr)',
\quad \mathfrak{Z}(t):={t^{k/\alpha}}{|\hat h(t)|},$$
so that if $h'(t)\ge0$, then it turns out that $\mathfrak{Z}$  is non increasing, being $\hat h(t)<0$ in $\mathbb R^+$. Thus, $0<\mathfrak{Z}(t) \le \mathfrak{Z}(0)=0$ in $\mathbb R^+$, a contradiction. So, if $\hat{h}(T_M)<0$, then $h'(t) <0$ in $\mathbb R^+$.

 In conclusion, there is always a right neighborhood of $0$, say $(0,P_0)$, with $P_0\in(0,\infty)$ opportune, in which $h$ is negative. 

Now, consider a cutoff function $\tau\in C^{\infty}(\mathbb R_0^+)$, nonincreasing and such that
$$\tau(t)=1, \,\, \text{if}\,\,0\le t\le P_{0} \quad \text{and} \quad \tau (t)=0 \,\, \text{if}\,\, t\ge P_{1},$$
for some real positive number $P_1>P_0$ and $P_0$ given above.
For instance, in the case in which $\hat{h}(T_M)> 0$,
we can choose $P_0=T_0$ and $P_1=T_1$, respectively the first and the second zero of $h$, whose existence is guaranteed by the behavior of $h$, cfr. Figure \ref{h1}. Differently, if $\hat{h}(T_M)=0$, then 
$h(T_M)=h'(T_M)=0$, and we can choose 
$P_0=T_M$, while $P_1$ any point in $(T_M,\infty)$, cfr. Figure~\ref{h2}.

Finally,
if $\hat{h}(T_M)< 0$ we can choose $P_0$ and $P_1$ in $\mathbb R^+$ without any restriction, cfr. Figure~\ref{h3}. 
We point out that, by \eqref{hT}, this latter case occurs roughly either if $\lambda>0$ arbitrary and $\beta$ large or if $\beta>0$ arbitrary and $\lambda$ large. Precisely,
define $\lambda^*_T$ and $\beta^*_T$ as follows
\begin{equation}\label{ltt}
\lambda^*_T= k(p^*-p)\,\biggl(\frac{p^*(\alpha p-k)}{\beta \|K\|_\infty } \biggr)^{(\alpha p-k)/\alpha(p^*-p)}\biggl(\frac{S\alpha }{2p(\alpha p^*-k) }\biggr)^{(\alpha p^*-k)/\alpha(p^*-p)}\cdot \frac 1{\|V\|_r}
\end{equation}
\begin{equation}\label{btt}
\beta^*_T= p^*(\alpha p-k)
\biggl(\frac {k(p^*-p)}{\lambda \|V\|_r} \biggr)^{\alpha(p^*-p)/(\alpha p-k)}\biggl(\frac{S\alpha }{2p(\alpha p^*-k)}\biggr)^{(\alpha p^*-k)/(\alpha p-k)}\cdot\frac 1{\|K\|_\infty}.
\end{equation}
Consequently, $\hat{h}(T_M)< 0$ holds either for all $\lambda>0$ and $\beta>\beta^*_T$  or  for all $\beta>0$ and $\lambda>\lambda^*_T$.
This property will appear crucial in the proof of Theorem \ref{1.1}.

\begin{figure}
     \centering
     \begin{subfigure}[b]{0.2\textwidth}
         \centering
         \includegraphics[width=\textwidth]{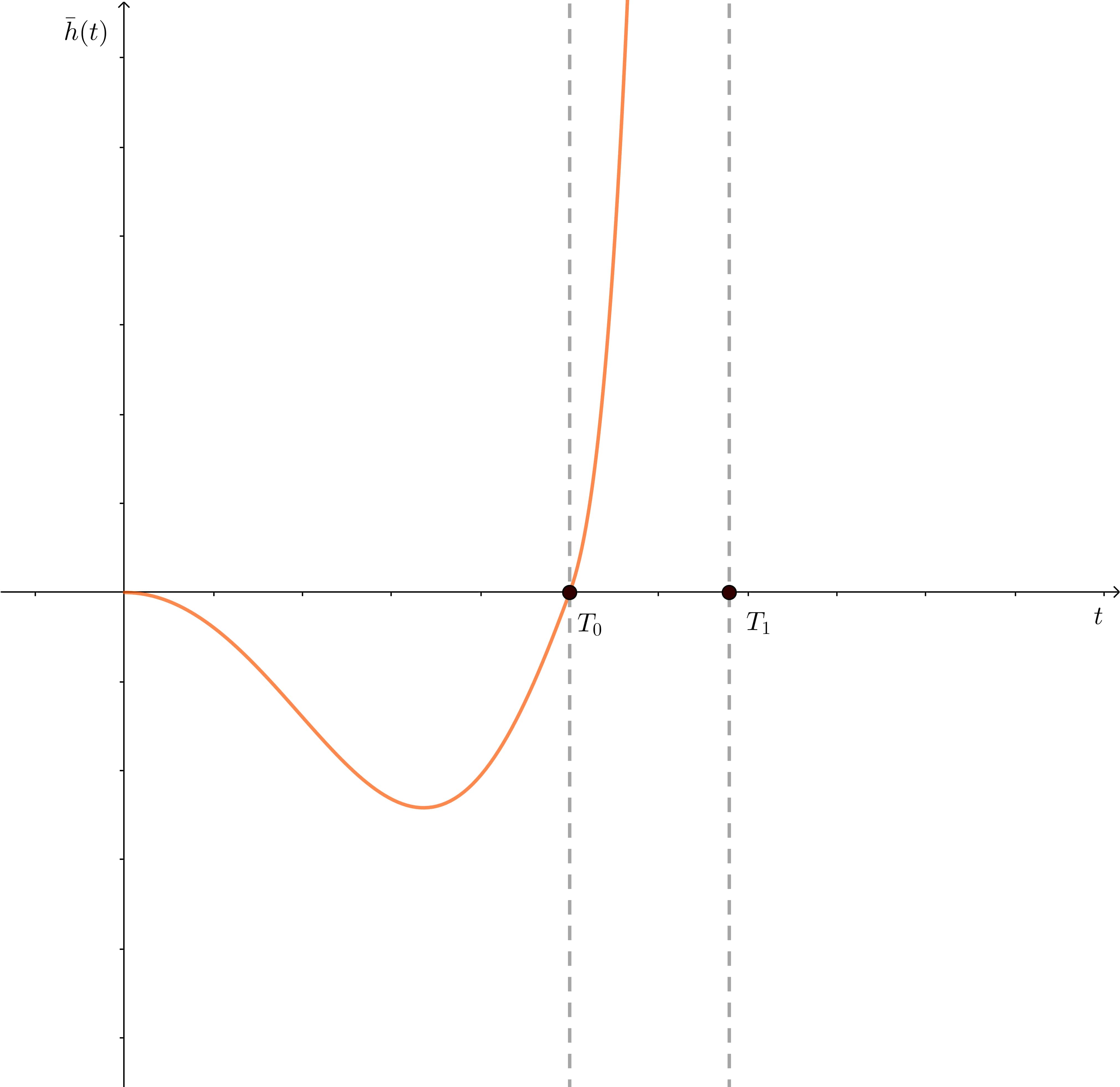}
         \caption{$\bar h$ if $\hat{h}(T_M)> 0$ }
         \label{h1t}
     \end{subfigure}
     \hfill
     \begin{subfigure}[b]{0.2\textwidth}
         \centering
         \includegraphics[width=\textwidth]{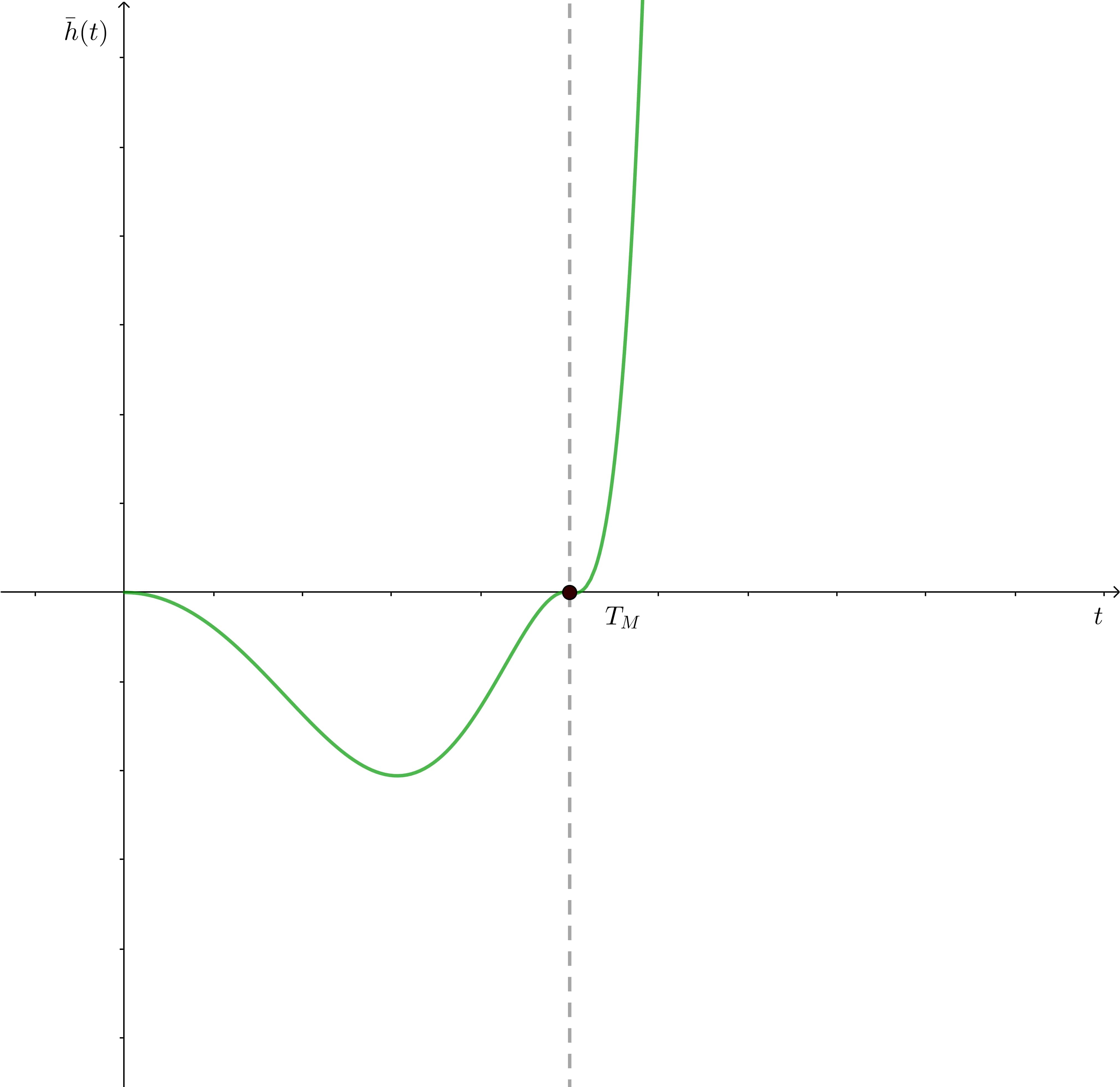}
         \caption{$\bar h$ if $\hat{h}(T_M)= 0$ }
         \label{h2t}
     \end{subfigure}
     \hfill
     \begin{subfigure}[b]{0.2\textwidth}
         \centering
         \includegraphics[width=\textwidth]{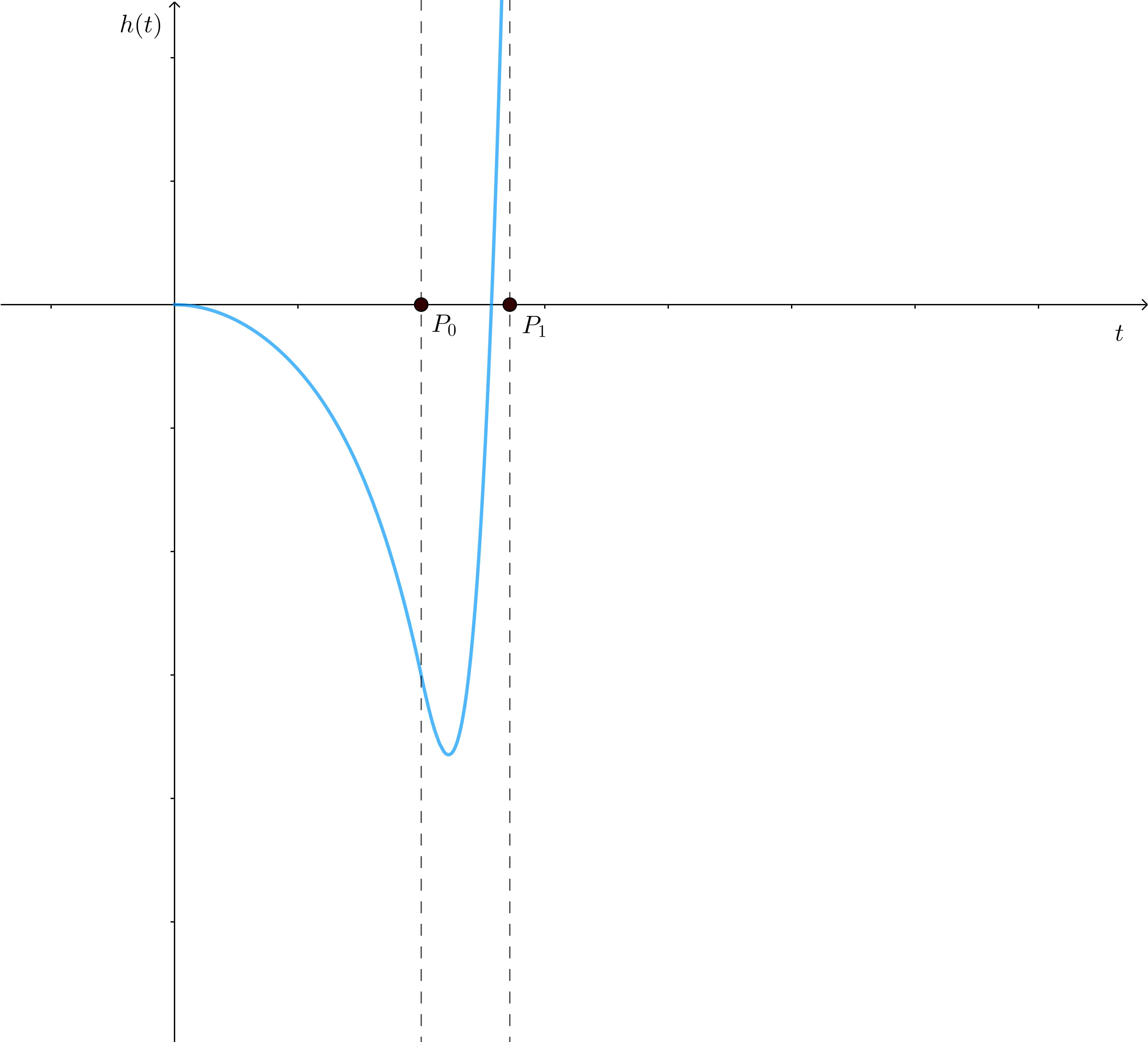}
         \caption{$\bar h$ if $\hat{h}(T_M)< 0$ }
         \label{h3t}
     \end{subfigure}
        \caption{$\bar h(t)$}
        \label{fig}
\end{figure}

By virtue of the cut-off function $\tau$ it is possible to define the truncated functional 
$$F_{\infty}(v)=\frac{1}{p}\|D v\|_p^{p}-\frac{\lambda}{k}\int_{\mathbb{R}^N}V|G^{-1}(v)|^{k}dx
-\beta\frac{\tau \left(\|v\|_{D^{1,p}}\right)}{\alpha p^{*}} \int_{\mathbb{R}^N}K|G^{-1}(v)|^{\alpha p^{*}}dx.$$
As above, we associate the real function
$\overline{h}(t)=\frac 1pt^{p}-\lambda c_Vt^{k/\alpha}-\beta c_Kt^{p^{*}}\tau (t)$, $t\in\mathbb R_0^+,$
whose behaviour, can be represented for instance by cfr. Figure \ref{h1t} if $P_0=T_0$, $P_1=T_1$ and by Figures \ref{h2t}, \ref{h3t} for the corresponding cases to Figures \ref{h2}, \ref{h3}.

\noindent
In particular, we have $F_{\infty}(v)\ge \overline{h}\left(\|v\|_{D^{1,p}}\right)$ for all $v\in D^{1,p}(\mathbb R^N)$ and
$F_{\lambda}(v)=F_{\infty}(v)$ if $0\le \|v \|_{D^{1,p}}\le P_{0}$.
Moreover, by the regularity both of the cut-off $\tau$ and of $F_\lambda$, we get $F_{\infty}\in C^{1}(D^{1,p}(\mathbb R^N),\mathbb{R})$ satisfying (PS)$_{c}$  for $\lambda<\lambda^*$ or $\beta<\beta^*$, as stated below, with $$\lambda^*=\min\{\lambda^*_{PS}, \lambda^*_T\}=\begin{cases}\lambda^*_{PS}
&\quad\text{ if }\,\,\biggl (\dfrac{\alpha p-k}p\biggr)^{p }>\biggl(\dfrac {\alpha p^*-k }{\alpha p^*}\biggr)^{p^*}\\
\lambda^*_{T}
&\quad\text{ if }\,\,\biggl (\dfrac{\alpha p-k}p\biggr)^{p }<\biggl(\dfrac {\alpha p^*-k }{\alpha p^*}\biggr)^{p^*}
\end{cases}$$ and
$$\beta^*=\min\{\beta^*_{PS}, \beta^*_T\}
=\begin{cases}\beta^*_{PS}
&\quad\text{ if }\,\,\biggl (\dfrac{\alpha p-k}p\biggr)^{p }>\biggl(\dfrac {\alpha p^*-k }{\alpha p^*}\biggr)^{p^*}\\
\beta^*_{T}
&\quad\text{ if }\,\,\biggl (\dfrac{\alpha p-k}p\biggr)^{p }<\biggl(\dfrac {\alpha p^*-k }{\alpha p^*}\biggr)^{p^*}
\end{cases}$$
with $\lambda^*_{PS}$, $\beta^*_{PS}$ given in \eqref{lambdastar_up}, \eqref{betastar_up}, while $\lambda^*_{T}$, $\beta^*_{T}$ defined in \eqref{ltt}, \eqref{btt}. 
Choosing either $\lambda$ or $\beta$ small is equivalent to fall in the cases described in Figures \ref{h1}-\ref{h1t} and \ref{h2}-\ref{h2t}.

\begin{lemma}\label{lem8}
Let $F_{\infty}$ be the truncated functional of $F_\lambda$. 
\begin{enumerate}
\item[(a)] For all $\lambda>0$ and $\beta\in (0,\beta^{*})$, then $F_{\infty}$ satisfies the
(PS)$_{c}$ with $c<0$.
\item[(b)] For all $\beta>0$ and $\lambda\in (0, \lambda^{*})$, then $F_{\infty}$ satisfies the
(PS)$_{c}$ with $c<0$.
\end{enumerate}
\end{lemma}

\begin{proof}
Assume either $\lambda>0$ and $\beta\in (0,\beta^{*})$ or $\beta>0$ and $\lambda\in (0, \lambda^{*})$, then it holds
\begin{enumerate}
\item[$(\mathfrak h)$] If $F_{\infty}(v)<0$, then $\|v\|_{D^{1,p}}<P_{0}$ and $F_{\lambda}(u)=F_{\infty}(u)$ for all $u$ in a
small enough neighbourhood of $v$.
\end{enumerate}
Indeed, choosing $\lambda$, $\beta$ as above, it immediately follows that $\hat h(T_M)\ge0$, so that the existence of $P_0$, with $h(P_0)=0$, is guaranteed, cfr. Figures \ref{h1}, \ref{h2}.
Consequently, we can apply the proof of Lemma 9 in \cite{BBF}, provided that  the change of variables  is taken into account.
\end{proof}

\section{Proof of Theorem \ref{1.1}}\label{main1}

This section is devoted to the proof of 
Theorem \ref{1.1}.
Before stating the symmetric Mountain Pass Theorem, we need the notion of the genus. 
\begin{definition}\label{defgen}
Let 
$\Sigma=\left\{ A\subset Y \backslash\left\{0\right\} | \,\ A \ \text{closed and symmetric}\
u\in A\Rightarrow -u\in A\right\}$, with $Y$ a real Banach space.
Let $A\in\Sigma$, the genus $\gamma\left( A\right)$ of $A$ is defined as the smallest integer $N$ such that there exists
$\Phi\in C\left( Y,\mathbb{R}^N\backslash\left\{0\right\}\right)$ such that $\Phi$ is odd and $\Phi\left(x\right)\ne 0$ for all
$x\in A$. We set $\gamma\left(\emptyset \right)=0$ and $\gamma\left(A\right)=\infty$ if there are no integers
 with the above property. 
\end{definition}

For an even functional on a Banach space, the symmetric Mountain Pass Theorem due to Ambrosetti–Rabinowitz \cite{AR}, Rabinowitz \cite{R} and Clark \cite{C}, as follows.

\begin{theorem} \label{t0}
Let $(Y, \|\cdot\|_Y)$ be an infinite dimensional Banach space and $E\in C^1(Y,\mathbb R)$ satisfying the following assumptions.
\begin{itemize}
    \item[{\bf (C1)}] $E$ is even, bounded from below, $E(0)=0$ and $E$ satisfies the Palais–Smale condition.
   \item[{\bf (C2)}] For each $j\in \mathbb N$, there exists an $A_j\in \Sigma_j$ such that $\sup_{u\in A_j} E(u)<0$.
\end{itemize}
Define 
$c_{j}=\inf_{A\in\Sigma_{j}}\sup_{u\in A}E(u)$,
$\Sigma_{j}=\left\{A\subset Y\backslash\left\{0\right\}: A \ \text{is closed in} \ Y,\,\, -A=A,\,\,\gamma(A)\ge j\right\}$.

Then each $c_j$ is a critical value of $E$ with $c_j\le c_{j+1} < 0$ for $j\in \mathbb N$ and $(c_j)_j$ converges to zero. 
Moreover, if $c_j=c_{j+1}=\dots=c_{j+i}=c$, then $\gamma(K_c)\ge i+1$, where $K_c$ is defined by $K_c=\{u\in Y : E'(u)=0, E(u)=c\}$.
\end{theorem}

Thus, for an even functional on an infinite dimensional Banach space, the symmetric Mountain Pass Theorem gives a
sequence of critical values which converges to zero. Under the same assumptions on the
functional, Kajikiya in \cite{K5} establishes the following critical point theorem which provides a sequence of critical points converging to zero. 

\begin{theorem}[Theorem 1, \cite{K5}] \label{t11}
Let $(Y, \|\cdot\|_Y)$ be an infinite dimensional Banach space and $E\in C^1(Y,\mathbb R)$ satisfying
 {\bf (C1)} and {\bf (C2)}, then one of the following holds.
\begin{itemize}
    \item There exists a sequence $(u_k)_k$ such that $E'(u_k)=0$, $E(u_k)<0$ and $\displaystyle{\lim_{k\to\infty} \|u_k\|_Y=0}$.
\item There exist two sequences $(u_k)_k$ and $(v_k)_k$ such that:
\begin{enumerate}
\item [(i)]
$E'(u_k)=0$, $E(u_k)=0$, $u_k\neq0$, $\displaystyle{\lim_{k\to\infty} \|u_k\|_Y=0}$
\item [(ii)] $E'(v_k)=0$, $E(v_k)<0$, $\displaystyle{\lim_{k\to\infty} E(v_k)=0}$ and  $\displaystyle{\lim_{k\to\infty} \|v_k\|_Y=c}$, with $c\neq 0$.
\end{enumerate}
\end{itemize}
\end{theorem}


\smallskip 
{\it Proof of Theorem \ref{1.1}.} 
In order to get our multiplicity result, we need to prove conditions {\bf (C1)} and {\bf (C2)} of Theorem \ref{t0}, as well as Theorem \ref{t11}, for the functional $F_\infty$. In particular, $F_\infty$ is even, bounded from below and satisfies the (PS)$_c$ condition for all $c<0$, see Section \ref{truncated}, thus {\bf (C1)} is proved. In order to prove {\bf (C2)} consider $K_c=K_{c, F_{\infty}}=\left\{u\in D^{1,p}(\mathbb R^N):\, F_{\infty}(u)=c, \,\,F'_{\infty}(u)=0\right\}\,$ defined in Theorem \ref{t0}
and take $m\in\mathbb N^+$. For $1\le j\le m$  let $c_j$ and $\Sigma_{j}$ as in Theorem \ref{t0}.
Our claim consists in proving that
\begin{equation}\label{cj<0}
c_j<0\quad\text{for all } j\ge1.
\end{equation}
To reach \eqref{cj<0} it is enough to prove that for all $j\in \mathbb{N}$, there exists an
$\varepsilon_{j}=\varepsilon(j)>0$ s.t.
\begin{equation}\label{gammaeinfty}
\gamma(F_{\infty}^{-\varepsilon_{j}})\ge j, \,\ \textnormal{where} \,\ F_{\infty}^{a}=\left\{v\in D^{1,p}(\mathbb R^N):\,
F_{\infty}(v)\le a\right\}\quad \text{with} \quad a\in\mathbb R.
\end{equation}
Let $\Omega_V\subset\mathbb{R}^N$, $|\Omega_V|>0$,  be a bounded open set
where  $V>0$ and extend functions in $D^{1,p}_{0}(\Omega_V)$  by
$0$ outside $\Omega_V$, where $D_{0}^{1,p}(\Omega_V)$ is the closure of $C^{\infty}_{0}(\Omega_V)$
in the norm $\|v\|_{D^{1,p}_{0}(\Omega_V)}$.
Take $W_{j}$ a $j$-dimensional subspace of $D^{1,p}_{0}(\Omega_V)$, thus all the norms in $W_{j}$ are equivalent.
For every $v\in W_{j}\setminus\{0\}$, we write $v=r_j w$ with $w\in W_{j}$ and $\|w\|_{D^{1,p}_{0}(\Omega_V)}=1$. In what follows we specify the choice of $r_j$. 
Note that there exists $d_{j}>0$ such that
\begin{equation}\label{dj} 
d_{j}\le 
\int_{\Omega_V}V|w|^{k}dx \end{equation}
 In particular, $V|w|^{k}\in L^1(\Omega_V)$
being $\Omega_V$ bounded, $V\in C(\mathbb R^N)$ and $k<p<p^*$.

By Lemma \ref{gprop}-{\bf a)}, for $\varepsilon>0$ there exists $\sigma=\sigma(\varepsilon)$ sufficiently small such that 
\begin{equation}\label{propg123}
|G^{-1}(t)|\ge \frac{1}{2^{1/k}}|t| 
\quad\text{if } |t|\le \sigma.
\end{equation}
Thus, since $\alpha<k<p$, by \eqref{dj} and \eqref{propg123}, we get
$$\begin{aligned}\int_{\Omega_V}V|G^{-1}(v)|^{k}dx&=\int_{\Omega_V}V\biggl( \frac12|v|^{k}+|G^{-1}(v)|^{k}-\frac12|v|^{k}\biggr)dx
\\&\ge \frac12\int_{\Omega_V}V|v|^{k}dx
+\int_{\{|v|\ge \sigma\}\cap\Omega_V}V\biggl(|G^{-1}(v)|^{k}
-\frac12|v|^{k}\biggr)dx\\
&\ge \frac{1 }{2}d_jr_j^{k}-\frac12\int_{\{|v|\ge \sigma\}\cap\Omega_V}V|v|^{k-p^*+p^*}dx
\\&\ge \frac12 d_jr_j^{k}-\frac12\|V\|_{L^\infty(\Omega_V)}
\sigma^{k-p^*}\|v\|_{p^*}^{p^*}=\frac 12 d_jr_j^{k}-cr_j^{p^*}, \quad c>0.
\end{aligned}$$
Consequently, by Lemma \ref{gprop}-{\bf g)}, the following holds 
$$\begin{aligned}
F_{\lambda}(v)&=\frac{1}{p}\int_{\Omega_V} |Dv|^p dx-\frac{\lambda }{k}
\int_{\Omega_V}V|G^{-1}(v)|^{k}dx-\frac{\beta}{\alpha p^{*}}\int_{\Omega_V}K|G^{-1}(v)|^{\alpha p^{*}}dx\\
&\le \frac{1}{p} r_j^p-\frac \lambda k\biggl (\frac 12 d_jr_j^{k}-cr_j^{p^*}\biggr)
+\frac{\beta 2^{(p-1)p^*}}{\alpha p^{*}}\|K\|_\infty\int_{\Omega_V}|v|^{p^{*}}dx
\\ &\le 
r_j^{k}\biggl[\frac{1}{p}r_j^{p-k}-\frac{\lambda}{2k}  d_j+Cr_j^{p^*-k}\biggr],
\end{aligned}$$
where $C>0$. Consequently,  for every $v\in W_{j}$, $v\neq 0$, we can choose $r_j\in(0,P_0)$ sufficiently small such that $\|v\|=r_j$ and, since $\alpha<k<p$, we obtain
$F_{\infty}(v)=F_\lambda(v)\leq -\varepsilon_j<0$,
by virtue of Lemma \ref{lem8}. Now, letting $S_{r_j}=\left\{v \in D^{1,p}(\mathbb R^N): \|v\|_{D^{1,p}_{0}(\Omega_V)}=r_j\right\},$ then $S_{r_j}\cap W_{j}\subset
F_{\infty}^{-\varepsilon_{j}}$. By Definition \ref{defgen}, $\gamma(F_{\infty}^{-\varepsilon_{j}})\ge \gamma(S_{r_j}\cap W_{j})=
\gamma(S^{j-1})=j$,
which proves claim \eqref{gammaeinfty}. Thus, from $F_{\infty}^{-\varepsilon_{j}}\in\Sigma_{j}$, we obtain
$c_j\le \sup_{u\in F_{\infty}^{-\varepsilon_{j}}}F_{\infty}(v)\le -\varepsilon_{j}<0.$

In particular, applying Theorem \ref{t11}, we get a sequence of solutions with negative energy of problem \eqref{Psc} such that, in both cases,  either $\|u_n\|\to0$  or $\|u_n\|\to c\neq0$ and $E_\lambda(u_n)\to0$ as $n\to\infty$.
In addition, in the latter case, there is also another sequence of nontrivial solutions $(\bar u_n)_n$ such that $E_\lambda(\bar u_n)=0$ and $\|\bar u_n\|\to0$ as $n\to\infty$.
This concludes the proof.
\qed

\begin{remark}\label{k_restr}
The case $p<k<\alpha p$ is not covered in Theorem \ref{1.1} because of some considerable difficulties. Indeed, to prove \eqref{cj<0} we need that $F_\infty$ has to be negative near $0$. Unfortunately, if $p<k$ it holds the opposite condition. Indeed, for any $v=r_j w\in W_j$ with $\|w\|_{D_0^{1,p}(\Omega_V)}=1$ and $r_j$ positive suitable, using Lemma \ref{gprop}-{\bf c), g)} and H\"older's inequality with exponents $p^*/k$ and $p^*/(p^*-k)$, we have
$$\begin{aligned}
F_\lambda(v)&\ge \frac{1}{p} r_j^p -\frac{\lambda}{k}\int_{\Omega_V} V|v|^k dx -\frac{\beta 2^{(p-1)p^*}}{\alpha p^*}\int_{\Omega_V} K|v|^{p^*}dx\\
&\ge \frac{1}{p} r_j^p -\frac{\lambda}{k}S^{-k/p}\|V\|_{L^{\infty}(\Omega_V)}|\Omega_V|^{(p^*-k)/p^*}r_j^k  -\frac{\beta 2^{(p-1)p^*}}{\alpha p^*}S^{-p^*/p}\|K\|_{L^\infty(\Omega_V)}r_j^{p^*},
\end{aligned}$$
since $V, K\in C(\mathbb R^N)$ and $\Omega_V$ is bounded, see the proof of Theorem \ref{1.1}.
Consequently, if $r_j$ is sufficiently small, say $r_j<P_0$, and $k>p$, then, by definition, $F_\infty(v)=F_\lambda(v)>0$.

On the other hand,  as it will be clear from  the proof of Theorem \ref{T_E>0} below, assuming $K\ge 0$ in $\mathbb R^N$, then one can obtain
$F_\lambda(v)<0$  for $r_j$ sufficiently large, since  $p, k/\alpha<p^*$. This choice of $r_j$ fits with the behaviour  of the functions $h,\bar h$ given in Figure \ref{h3} and \ref{h3t}, see Section $\ref{truncated}$ for details, namely when
$h$ is nonincreasing.
However, this latter case cannot occur since 
(PS)$_c$ property for $F_\lambda$ follows from Lemma \ref{lem6}
by requiring either $\lambda$ fixed and $\beta$ small or viceversa, in turn $\hat h(T_M)\ge 0$ by \eqref{hT}, that are situations described in Figures \ref{h1}-\ref{h1t} and \ref{h2}-\ref{h2t}.
 
\end{remark}

\section{Proof of Theorem  \ref{T_E>0} }\label{main2}

This section is devoted to the proof of the result under a symmetric setting. In particular, Theorem \ref{T_E>0}, whose statement is given in the Introduction, gives infinitely many solutions for problem \eqref{PS} with no assumptions on the parameters $\lambda$ and $\beta$. Indeed, thanks to the symmetric environment considered, by the corollary below, the (PS)$_c$ condition for the energy functional $F_\lambda$ holds in a very general setting.

\begin{corollary}\label{corpsc}
Let $\alpha<k<\alpha p^*$.
If \eqref{K0Kinfty=0} holds, then the functional $F_\lambda$ satisfies (PS)$_c$ condition in $D^{1,p}_T(\mathbb R^N)$ for every $c\in \mathbb R$.
\end{corollary}
For the proof of the result above we refer to Corollary 1 in \cite{BBFS}.

\begin{remark}\label{no_mult}
Unfortunately, we cannot improve Theorem \ref{1.1} in a symmetric setting in view of the new property given by Corollary  \ref{corpsc}, since Lemma \ref{lem8} seems difficult to be extended 
to any $\lambda, \beta>0$.
Indeed, in the proof of Lemma \ref{lem8} we deduce the validity of (PS)$_c$ condition for the truncated energy functional $E_\infty$ from property $(\mathfrak h)$. This latter property is easily obtained for the cases \ref{h1t}, \ref{h2t}, valid roughly for either $\lambda$ or $\beta$ small, but, for the case in Figure \ref{h3t} it fails since, if $F_\infty(v)<0$, it could be possible that $\|v\|_{D^{1,p}}>P_0$, yielding $F_\infty(v)\neq F_\lambda(v)$.
\end{remark}

Now we are ready to prove Theorem \ref{T_E>0}, by using the Fountain Theorem.

{\it Proof of Theorem \ref{T_E>0}.}
The proof consists in applying Theorem \ref{FT} with $\mathcal{G}=\mathbb{Z}/2$ and $M=D^{1,p}_T(\mathbb{R^N})$.
Note that, by Remark \ref{aft}, assumption {\bf (A1)} is verified.
The energy functional $F_{\lambda}$, from Corollary \ref{corpsc}, satisfies the (PS)$_{c}$ condition for every $c\in\mathbb{R}$, so that also {\bf (A4)} of Theorem \ref{FT} holds.
Since $0\not\equiv K\ge0$ in $\mathbb R^N$ and $K\in C(\mathbb R^N)$, there exists an open bounded subset $\Omega_K$ of $\mathbb R^N$ with
$K>0$ in $\Omega_K$. By the $T$-symmetry of $K$, then also $\Omega_K$ is $T$-symmetric, in turn we can define $D^{1,p}_T(\Omega_K)$. Now,   extending functions in $D^{1,p}_T(\Omega_K)$ by $0$ outside $\Omega_K$, we can consider  $D^{1,p}_T(\Omega_K)\subset D^{1,p}_T(\mathbb R^N)$. 
Assume $(Y_{m})_m$ be an increasing sequence of subspaces of $D^{1,p}_T(\Omega_K)$ with $dim(Y_m)=m$. 
In particular, if $v\in Y_{m}$, $v\ne 0$, then we can write $v=\rho_{m}\omega$ with $\omega\in Y_m$ such that $\|\omega\|=1$, so that $\rho_m=\|v\|$.
Moreover, there exists a constant
$\varepsilon_{m}>0$ such that 
\begin{equation}\label{4.4.1}
\int_{\mathbb R^N}K|\omega|^{p^{*}}dx=\int_{\Omega_K}K|\omega|^{p^{*}}dx\ge\varepsilon_{m}.
\end{equation}
By Lemma \ref{gprop}-{\bf b)}, for all $\varepsilon>0$ there exists $M>0$ large enough, such that $|G^{-1}(t)|^\alpha\ge p^{-1/p^*} |t|$ for $|t|\ge M$ choosing $\varepsilon=2^{1/p}-p^{-1/p^*}$, so that, thanks to  Lemma \ref{gprop}-{\bf g)}, \eqref{4.4.1} and since $K, V\ge 0$, the following holds 
$$\begin{aligned}
F&_{\lambda}(v)=
\frac{1}{p}\rho_{m}^{p}-\frac{\lambda}{k}\int_{\Omega_K} V |G^{-1}(v)|^k dx -\frac{\beta}{\alpha p^*}\int_{\Omega_K} K\biggl[ \frac{|v|^{p^*}}{p} +|G^{-1}(v)|^{\alpha p^*}-\frac{|v|^{p^*}}{p}\biggr] dx \\
 &=\frac{1}{p}\rho_{m}^{p}-\frac{\lambda}{k}\int_{\Omega_K} V |G^{-1}(v)|^k dx -\frac{\beta}{\alpha pp^*}\int_{\Omega_K} K|v|^{p^*}dx
 \\&\quad-\frac{\beta}{\alpha p^*}\int_{\{|v|\ge M\}\cap\Omega_K}K \biggl[|G^{-1}(v)|^{\alpha p^*}-\frac{|v|^{p^*}}{p} \biggr] dx  \\
 &\quad- \frac{\beta}{\alpha p^*}\int_{\{|v|< M\}\cap\Omega_K} K|G^{-1}(v)|^{\alpha p^*}dx+ \frac{\beta}{\alpha pp^*}\int_{\{|v|< M\}\cap\Omega_K} K |v|^{p^*} dx\\
 &\le \frac{1}{p}\rho_{m}^{p} -\frac{\beta}{\alpha p p^*}\int_{\Omega_K} K |v|^{p^*} dx + \frac{\beta}{\alpha p p^*}\int_{\{|v|< M\}\cap\Omega_K} K |v|^{p^*} dx \\
&\le \frac{1}{p}\rho_{m}^{p}
-\frac{\beta}{\alpha p p^*}
\int_{\Omega_K} K |v|^{p^*} dx + \frac{\beta}{\alpha p p^*}\|K\|_\infty M^{p^*-k/\alpha}\int_{\Omega_K}|v|^{k/\alpha} dx\\
&\le \frac{1}{p}\rho_{m}^{p}-
\frac{\beta\varepsilon_{m}}{\alpha p p^{*}}\rho_{m}^{p^{*}}+C\rho_{m}^{k/\alpha}< 0 \end{aligned}$$
for sufficiently large $\rho_{m}$, since $k<\alpha p^*$, where
$C=\frac{\beta}{\alpha p p^*}\|K\|_\infty M^{p^*-k/\alpha}|\Omega_K|^{1-k/\alpha p^*}S^{-k/\alpha p}.$
This proves {\bf (A2)} of Theorem \ref{FT}.
Condition {\bf (A3)} follows exactly as in \cite{BBFS}, taking into account the properties of $v$. 
Then applying Theorem \ref{FT}, the energy functional $F_{\lambda}$ has unbounded  sequence of critical values in $D^{1,p}_T(\mathbb{R^N})$, so that
Theorem \ref{T_E>0} is proved.
\qed

\section{Further results in the singular case: $0<\alpha<1$}\label{sing}

In this section we extend Theorem 1.1 in \cite{BFsc} where  the so called singular version of equation \eqref{Psc}, i.e. $0<\alpha<1$, is considered, and for which the critical exponent is $p^*$,
\begin{equation}\label{Pscc}
-\Delta_{p}u-\frac{\alpha}{2}\Delta_p(|u|^\alpha)|u|^{\alpha-2}u=\lambda V(x)|u|^{k-2}u
+\beta K(x)|u|^{p^{*}-2}u \quad\text{in }\mathbb{R}^N,
\end{equation}
where $1<k<p^*$. 
This case is known as "singular" since the function $g(t)$ in \eqref{gpicc} is singular for $t=0$, differently from the case $\alpha>1$ where $g$ is coercive at $\infty$. 

The natural energy functional associated with equation \eqref{Pscc} is 
$$\hat E_{\lambda}(v)=\frac{1}{p}\int_{\mathbb R^N}g(u)^p|D u|^{p} dx
-\frac{\lambda}{k}\int_{\mathbb{R}^N}V|u|^{k}dx-\frac{\beta}{p^{*}}\int_{\mathbb{R}^N}K|u|^{p^{*}}d,$$
while, the functional we have to deal with, after the same change of variable described in Section \ref{prel3}, is the following
$$\hat F_{\lambda}(v)=\frac{1}{p}\int_{\mathbb R^N}|D v|^{p} dx-\frac{\lambda}{k}
\int_{\mathbb{R}^N}V|G^{-1}(v)|^{k}dx-\frac{\beta}{p^{*}}\int_{\mathbb{R}^N}K|G^{-1}(v)|^{p^{*}}dx.$$
To manage the singular case, we have to take care of the different growth of $G^{-1}$ at $0$ and at $\infty$ respect to the case $\alpha>1$.
Multiplicity results for solutions of \eqref{Pscc} with negative energy in the singular case are obtained using the same technique as for Theorem \ref{1.1} but under the more restrictive condition  $2<k<\alpha p$. Here, we enlarge the interval for $k, p$ up to consider $1<k<\alpha p$ and we add properties on the behaviour of solutions, as follows.

\begin{theorem}\label{1.1s}
Let $N\ge3$ and $1<k<\alpha p<p< N$. Assume that $K$ satisfies \eqref{K(x)1} and $0\le V\in  L^{\varrho}(\mathbb{R}^{N})\cap C(\mathbb R^N)$ with $\varrho=p^{*}/(p^{*}-k)$. Then, 
\begin{itemize}
\item[(i)] For any $\lambda>0$, there exists $\hat \beta>0$ such that for any $0<\beta<\hat\beta$, then equation \eqref{Pscc} has infinitely many nontrivial solutions $(u_n)_n\subset D^{1,p}(\mathbb R^N)$  such that $\hat E_\lambda(u_n)<0$ and $\|u_n\|\to0$ as $n\to\infty$.
\item[(ii)]For any $\beta>0$, there exists $\hat\lambda>0$ such that for any $0<\lambda<\hat\lambda$, then equation \eqref{Pscc} has infinitely many nontrivial solutions $(u_n)_n\subset D^{1,p}(\mathbb R^N)$  such that $\hat E_\lambda(u_n)<0$ and $\|u_n\|\to0$ as $n\to\infty$.
\end{itemize}
\end{theorem}

The smaller interval $2<k<\alpha p$ in \cite{BFsc} follows from the application of Lemma 3.3 in \cite{BFsc}, which is crucial in order to obtain the validity of (PS)$_c$ condition for $\hat F_{\lambda}$. We point out that the restriction $2<k<\alpha p$ is not required in proving that  concentration around points and at infinity cannot occur.
On the other hand, Lemma 3.3 in \cite{BFsc} cannot be applied in the case  $\alpha>1$  since the critical exponent becomes $\alpha p^*>p^*$ so different arguments are required. Surprisingly, these last tools, applied to the singular case, allow us to cover the entire interval $1<k<\alpha p$. 
In particular, it holds the following.

\begin{lemma}\label{lem6s}
Suppose $0\le V\in  L^{\varrho}(\mathbb{R}^{N})\cap C(\mathbb R^N)$ with $\varrho=p^{*}/(p^{*}-k)$ and $K$ satisfy \eqref{K(x)1}. Let $1<k<p$ and $c<0$. Then 
\begin{enumerate}
    \item[(I)] For any $\lambda>0$, there exists $\hat\beta>0$ defined as follows
$$\hat\beta=\frac{ \alpha}{\|K\|_\infty}\biggl(\frac{kp^*}{\lambda N \|V\|_r(p^*-k)}\biggr)^{p^2/(N-p)(p-k)} S^{(p^*-k)/(p-k)}$$
    such that for every $\beta\in(0,\hat\beta]$, then $\hat F_{\lambda}$ satisfies (PS)$_{c}$ condition.
    \item[(II)] For any $\beta>0$, there exists $\hat\lambda>0$ defined as follows
$$\hat\lambda=S^{(p^*-k)/(p^*-p)} \frac{kp^*}{N(p^*-k)}\cdot\frac{1}{\|V\|_r} \cdot \biggl(\frac{\alpha}{\beta\|K\|_\infty}\biggr)^{(p-k)/(p^*-p)},$$
    such that for every $\lambda\in(0,\hat\lambda]$, then $\hat F_{\lambda}$ satisfies (PS)$_{c}$ condition.
\end{enumerate}
\end{lemma}

\begin{proof}
Let $(v_n)_n$ be a (PS)$_c$ sequences for the functional $\hat F_{\lambda}(v)$. 
We now repeat word by word up the proof of Lemma 3.4 in \cite{BFsc} until formula (3.28) obtained thanks to the  application of the Brezis Lieb Lemma, precisely
\begin{equation}\label{g-11s}
\lim_{n\to\infty}\|G^{-1}(v_{n})-G^{-1}(v)\|_{p^{*}}=0.
\end{equation}
At this point, using that $\hat F_{\lambda}(v_n)\to \hat F_{\lambda}(v)=c$ by weak continuity of the functional, we have
$$\begin{aligned}
\int_{\mathbb R^N}&|D v_n|^{p} dx-\int_{\mathbb R^N}|D v|^{p} dx +o(1)\\&=
\frac{\lambda p}{k}\int_{\mathbb{R}^N}V\bigl[|G^{-1}(v_n)|^{k}-|G^{-1}(v)|^{k}\bigr]dx
+\frac{\beta p}{\alpha p^{*}}\int_{\mathbb{R}^N}K\bigl[|G^{-1}(v_n)|^{ p^{*}}
-|G^{-1}(v)|^{ p^{*}}\bigr]dx
\end{aligned}$$
as $n\to\infty$. Thus, by \eqref{g-11s} and H\"older's inequality, the right hand side of the above equation tends to $0$ as $n\to\infty$, so that 
 $\|Dv_n\|_p\to\|Dv\|_p$ as $n\to\infty$. Consequently, since $D^{1,p}(\mathbb R^N)$ is uniformly convex and by $Dv_n \rightharpoonup Dv$ in $\mathbb R^N$, by Proposition 3.32 in \cite{Brezis}, we get
the strong convergence in $L^p(\mathbb R^N)$ of the sequence $(D v_n)_n$.
\end{proof}


{\it Proof of Theorem \ref{1.1s}.}
It is enough to follow the proof Theorem 1 in \cite{BFsc} where, in place of Lemma 3.4, we use Lemma \ref{lem6s} above.
In addition, to obtain the asymptotic behaviours we refer  to Theorems \ref{t0}, \ref{t11} given in Section \ref{main1}. \hfill\qed

As described in Remark \ref{no_mult} it seems difficult to extend the above theorem for every $\lambda, \beta>0$ in the symmetric setting of Theorem \ref{T_E>0}, while for the corresponding multiplicity result for solutions with positive energy we refer to Theorem 1.3 in \cite{BFsc}.

\section*{Acknowledgments}
R. Filippucci and L. Baldelli are members of the {\em Gruppo Nazionale per
l'Analisi Ma\-te\-ma\-ti\-ca, la Probabilit\`a e le loro Applicazioni}
(GNAMPA) of the {\em Istituto Nazionale di Alta Matematica} (INdAM).
R. Filippucci was partly supported by  {\em Fondo Ricerca di Base di Ateneo Esercizio} 2017-19 of the University of Perugia, named {\em  Problemi con non linearit\`a dipendenti dal gradiente} and by INdAM-GNAMPA Project 2022 titled {\em Equazioni differenziali alle derivate parziali in fenomeni non lineari} (E55F22000270001).
L. Baldelli was partially supported by National Science Centre, Poland (Grant No. 2020/37/B/ST1/02742).

\end{document}